\documentclass{amsart}
\usepackage{latexsym,amsxtra,amscd,ifthen}
\usepackage{amsfonts}
\usepackage{verbatim}
\usepackage{amsmath}
\usepackage{amsthm}
\usepackage{amssymb}

\numberwithin{equation}{section}

\theoremstyle{plain}
\newtheorem{theorem}{Theorem}[section]
\newtheorem{lemma}[theorem]{Lemma}
\newtheorem{proposition}[theorem]{Proposition}
\newtheorem{corollary}[theorem]{Corollary}

\theoremstyle{definition}
\newtheorem{definition}[theorem]{Definition}
\newtheorem{example}[theorem]{Example}

\newtheorem{remark}[theorem]{Remark}
\newtheorem{question}[theorem]{Question}

\makeatletter              
\let\c@equation\c@theorem

\makeatother

\newcommand{\e}{\mathrm e}

\newcommand{\boplus}{\textstyle \bigoplus}
\newcommand{\bcup}{\textstyle \bigcup}

\DeclareMathOperator{\gldim}{gldim}

\DeclareMathOperator{\Ext}{Ext}
\DeclareMathOperator{\Tor}{Tor}
\DeclareMathOperator{\pdim}{projdim}

\DeclareMathOperator{\HB}{H}

\DeclareMathOperator{\cd}{cd}
\DeclareMathOperator{\lcd}{lcd}

\DeclareMathOperator{\gr}{gr}

\DeclareMathOperator{\htr}{Htr}
\DeclareMathOperator{\GKtr}{GKtr }
\DeclareMathOperator{\fhtr}{H_1tr}
\DeclareMathOperator{\shtr}{H_2tr}
\DeclareMathOperator{\thtr}{H_3tr}
\DeclareMathOperator{\htrgr}{Htr}
\DeclareMathOperator{\fhtrgr}{H_1tr}
\DeclareMathOperator{\shtrgr}{H_2tr}
\DeclareMathOperator{\thtrgr}{H_3tr}
\DeclareMathOperator{\Ktr}{Ktr}
\DeclareMathOperator{\tr}{tr}

\DeclareMathOperator{\ltr}{Ltr}

\DeclareMathOperator{\Kdim}{Kdim}
\DeclareMathOperator{\Cdim}{Cdim}
\DeclareMathOperator{\injdim}{injdim}

\DeclareMathOperator{\GKdim}{GKdim}
\DeclareMathOperator{\End}{End}
\DeclareMathOperator{\R}{R}
\DeclareMathOperator{\D}{{\sf D}}
\DeclareMathOperator{\Hom}{Hom}
\DeclareMathOperator{\RHom}{RHom}
\DeclareMathOperator{\GrMod}{{\sf GrMod}}

\DeclareMathOperator{\Mod}{{\sf Mod}}

\newcommand{\fm}{\mathfrak{m}}

\begin{document}

\title{Homological transcendence degree}

\author{Amnon Yekutieli and James J. Zhang}

\address{A. Yekutieli: Department of  Mathematics 
Ben Gurion University, Be'er Sheva 84105, Israel}

\email{amyekut@math.bgu.ac.il}

\address{J.J. Zhang: Department of Mathematics, 
Box 354350, University of Washington, Seattle, 
Washington 98195, USA}

\email{zhang@math.washington.edu}

\begin{abstract}
Let $D$ be a division algebra over a base field $k$. 
The homological transcendence degree of $D$, denoted 
by $\htr D$, is defined to be the injective dimension 
of the algebra $D\otimes_k D^{\circ}$. We show that 
$\htr$ has several useful  properties which the 
classical transcendence degree has. We extend some 
results of Resco, Rosenberg, Schofield and Stafford, 
and compute $\htr$ for several classes of division 
algebras. The main tool for the computation is 
Van den Bergh's rigid dualizing complex. 
\end{abstract}




\thanks{This research was supported by the US-Israel 
Binational Science Foundation. The second author 
was partially supported by the US National Science 
Foundation.}

\maketitle


\setcounter{section}{-1}
\section{Introduction}
\label{xxsec0}

Throughout $k$ is a commutative base field. By 
default all algebras and rings are $k$-algebras, 
and all homomorphisms are over $k$. This paper is 
mainly about division algebras that are infinite 
dimensional over their centers. Such division 
algebras appear naturally in noncommutative ring 
theory, and recently there are many new examples 
coming from noncommutative projective geometry. 
One important question in noncommutative 
algebra/geometry is the classification of division 
algebras of transcendence degree 2  (see some 
discussion in \cite{Ar}). Similarly to the 
commutative situation, the classification of 
division algebras of transcendence degree 2 would 
be equivalent to the birational classification of 
integral noncommutative projective surfaces. 
Quantum ${\mathbb P}^3$'s (the classification of 
which has not yet been achieved), will provide 
new examples of division algebras of transcendence 
degree 3. Other division algebras such as the 
quotient division rings of Artin-Schelter regular 
algebras will certainly play an important role in 
noncommutative projective geometry.

Most division algebras arising from noncommutative 
projective geometry should have finite transcendence 
degree. But what is the definition of transcendence 
degree for a division algebra infinite dimensional 
over its center? The first such definition is due 
to Gelfand and Kirillov \cite{GK}. Let $\GKdim$ denote 
the Gelfand-Kirillov dimension. Then the 
{\it Gelfand-Kirillov transcendence degree} of a 
division algebra $D$ is defined to be
$$\GKtr D=\sup_{V}\inf_{z} \; \GKdim k[z V],$$
where $V$ runs over all finite dimensional 
$k$-subspaces of $D$, and $z$ runs over all nonzero 
elements in $D$. This is probably the first simple 
invariant that distinguishes between the Weyl skew 
fields, since $\GKtr D_n=2n$, where $D_n$ is the 
$n$-th Weyl skew field. Note that $\GKdim D_n=\infty$ 
for all $n$, so it does not provide any useful 
information. Partly due to the complicated definition, 
$\GKtr$ is very mysterious. For example, it is not 
known whether $\GKtr D_1\leq \GKtr D_2$ when 
$D_1\subset D_2$  are division algebras. Also there 
are only a handful families of division algebras 
for which $\GKtr$ was computed explicitly 
\cite{GK, Lo, Zh1}. Recently the $\GKtr$ of the 
quotient division rings of twisted homogeneous 
coordinate rings was computed in 
\cite[Corollary 5.8]{RZ}. 

The second author gave another definition, 
called {\it lower transcendence degree}, denoted by 
$\ltr$ \cite{Zh2}. In general it is not clear whether 
$\ltr=\GKtr$. Several basic properties of classical 
transcendence degree have been established for $\ltr$. 
Using the properties of $\ltr$, one can compute both 
$\ltr$ and $\GKtr$ for several more classes of 
division algebras. For example, both $\ltr$ and 
$\GKtr$ of the quotient division ring of any 
$n$-dimensional Sklyanin algebra are equal to $n$. 
It is not hard to see from the definition that both 
$\GKtr$ and $\ltr$ are bounded  by $\GKdim$. So for 
some classes of rings one can obtain upper bounds for 
these two invariants. However it is not easy to compute 
the exact value in general. One open question is the 
following.

\begin{question} 
\label{xxque0.1}
Let $Q$ be the quotient division ring of a noetherian 
Artin-Schelter regular Ore domain of global dimension 
$n\geq 4$. Is $\ltr Q= \GKtr Q=n$?
\end{question}

The answer to Question \ref{xxque0.1} is ``yes'' for 
$n\leq 3$ (see \cite[Theorem 1.1(10)]{Zh1} and 
\cite[Theorem 0.5(4)]{Zh2}). Also it can be shown that 
both $\ltr Q$ and $\GKtr Q$ in Question \ref{xxque0.1} 
are finite.

It is fundamental to have a transcendence degree that 
is useful and computable for a large class of division 
algebras including those arising from noncommutative 
projective geometry. In this paper we introduce a new 
definition of transcendence degree which is defined 
homologically; and show that this transcendence 
degree is computable for many algebras with good 
homological properties, including all quotient 
division rings of Artin-Schelter regular algebras. 
Let $D$ be a division algebra over $k$. The 
{\it homological transcendence degree} of $D$ is 
defined to be 
$$\htr D=\injdim D\otimes_k D^\circ$$
where $D\otimes_k D^\circ$ is viewed as a left module 
over itself. Here $D^\circ$ is the opposite ring of 
$D$ and $\injdim$ denotes the injective dimension of 
a left module.

The idea of studying homological invariants of 
division algebras first appeared in Resco's papers 
\cite{Re1,Re3} for commutative fields. Later this 
was used by Stafford \cite{St} to study the quotient 
division ring of the group ring $kG$ of a torsion-free 
polycyclic-by-finite group $G$, and the quotient 
division ring of the universal enveloping algebra 
$U(\mathfrak g)$ of a finite dimensional Lie 
algebra $\mathfrak g$; and by Rosenberg \cite{Ro} 
to study the Weyl skew fields. Schofield extended 
this idea effectively to stratiform simple artinian 
rings, and proved several wonderful results in 
\cite{Sc1}. In addition there are other papers that 
studied various invariants of the tensor
product of divisions rings \cite{RSW,Sc3,Va}. Our 
definition of transcendence degree is motivated by 
the work of  Resco, Rosenberg, Schofield and Stafford. 
Using the results of Schofield \cite{Sc1} one can 
show that $\htr D=n$ where $D$ is a stratiform
simple artinian ring of length $n$ 
[Proposition \ref{xxprop1.8}]. A similar computation 
works for the division rings studied by Resco,
Rosenberg and Stafford. However we intend to cover 
a large class of division algebras for which the 
methods of the above people may not apply. 

If $A$ is an Ore domain, let $Q(A)$ denote the 
quotient division ring of $A$. We prove the following.

\begin{theorem}
\label{xxthm0.2} Let $A$ be an Artin-Schelter 
regular graded Ore domain. Then $\htr Q(A)=\gldim A$.
\end{theorem}

In Theorem \ref{xxthm0.2} above we do not assume that 
$A$ is noetherian. One of the main tools for this 
computation is the rigid dualizing complex, which we 
will review in Section \ref{xxsec3}. As a consequence, 
if $A$ and $B$ are two Artin-Schelter regular Ore 
domains of different global dimensions, then $Q(A)$ 
is not isomorphic to $Q(B)$. Theorem \ref{xxthm0.2} 
fails without the Artin-Schelter condition. By 
Propositions \ref{xxprop7.6} and \ref{xxprop7.8} 
below there is a connected graded Koszul Ore domain 
$A$ of GK-dimension 4 and global dimension 4 such that 
$$\htr Q(A)=3<4=\GKtr Q(A)=\GKtr A.$$

One big project would be to compute the homological transcendence 
degree for all division algebras that are not constructed from 
Artin-Schelter regular algebras. We start this task with the  
quotient division rings of connected graded domains with some 
mild homological hypotheses.

\begin{theorem}
\label{xxthm0.3} Let $A$ be a connected graded noetherian domain. 
Let $Q=Q(A)$.
\begin{enumerate}
\item 
If $A$ has enough normal elements, then $\htr Q=\GKdim A$.
\item
If $A$ has an Auslander balanced dualizing complex and if 
$A\otimes_k Q^{\circ}$ is noetherian, then $\htr Q=\cd A$.
Here $\cd$ denotes the cohomological dimension defined in
\textup{Definition \ref{xxdefn6.1}(d)}.
\item
If $A$ is an Artin-Schelter Gorenstein ring and if $A\otimes_k Q^{\circ}$ 
is noetherian, then $\htr Q=\injdim A$.
\end{enumerate}
\end{theorem}

In addition to Theorems \ref{xxthm0.2} and \ref{xxthm0.3} above 
we have Propositions \ref{xxprop5.4} and \ref{xxprop6.4} and 
Theorem \ref{xxthm6.9} which compute $\htr$ for filtered rings. 
These seemingly technical results cover several 
different classes of division rings. For example, $\htr$ is computable
for the following classes of division rings:

\begin{enumerate}
\item[(i)] 
quotient division rings of affine prime PI algebras, 
\item[(ii)] 
the Weyl skew fields,
\item[(iii)] 
quotient division rings of enveloping algebras  
$U(\mathfrak g)$ of finite dimensional Lie algebras, and 
\item[(iv)] 
quotient division rings of some other quantum algebras studied 
by Goodearl and Lenagan \cite{GL}. 
\end{enumerate}

Another advantage of this new definition is that it is easy to 
verify some useful properties similar to those of the classical 
transcendence degree. Let $\tr$ denote the classical transcendence 
degree of a commutative field (or a PI division algebra).

\begin{proposition}
\label{xxprop0.4} Let $D$ be a division algebra and $C$ be a 
division subalgebra of $D$.
\begin{enumerate}
\item 
$\htr C\leq \htr D$.
\item 
If $D$ is finite as a left (or a right) $C$-module, then 
$\htr C=\htr D$.
\item 
Suppose $C\otimes_k C^\circ$ is noetherian of finite global
dimension.
If $D$ is the quotient division ring of the skew polynomial ring 
$C[x;\alpha]$ for some automorphism $\alpha$ of $C$, 
then $\htr D =\htr C +1$.
\item 
If $D$ is a PI division ring and if the center is finitely generated 
over $k$ as a field, then $\htr D =\tr D$.
\item
$\htr D=\htr D^{\circ}$.
\end{enumerate}
\end{proposition}

As a consequence of Theorem \ref{xxthm0.2} and Proposition 
\ref{xxprop0.4}(a), if $A$ and $B$ are two Artin-Schelter
regular Ore domains and $\gldim A<\gldim B$, then 
there is no algebra homomorphism from $Q(B)$ to $Q(A)$.

\section{Definitions and basic properties}
\label{xxsec1}

Let $A$ be an algebra over the base field $k$. Let $A^{\circ}$ be the 
opposite ring of $A$, and let $A^\e$ be the enveloping algebra 
$A\otimes A^{\circ}$, where $\otimes$ denotes $\otimes_k$. Note that 
the switching operation $a\otimes b \mapsto b\otimes a$ extends to an 
anti-automorphism of the algebra $A^\e$. Usually we work with left 
modules. A right $A$-module is viewed as an $A^{\circ}$-module, and an 
$A$-bimodule is the same as an $A^\e$-module.
An $A$-module is called {\it finite} if it is finitely generated
over $A$. 

We have already seen the definition of $\htr$ in the introduction. 
To compute $\htr$ it is helpful to introduce a few related 
invariants, which are called modifications of $\htr$. For a ring
$B$ and a $B$-module $N$, we denote by $\injdim_B N$ (and $\pdim_B N$) 
the injective dimension (respectively, the projective dimension) of 
$N$ as a $B$-module. If $N=B$, we simplify $\injdim_B B$ to 
$\injdim B$. 

\begin{definition} 
\label{xxdefn1.1}
Let $A$ be a $k$-algebra. 
\begin{enumerate}
\item 
The {\it homological transcendence degree} of $A$ is defined to be
$$\htr A=\injdim  A^\e.$$
\item 
The {\it first modification} of $\htr$ is defined to be
$$\fhtr A=\sup\{i\;|\; \Ext^i_{A^\e}(A,A^\e)\neq 0\}.$$
\item 
The {\it second modification} of $\htr$ is defined to be
$$\shtr A=\sup\{\injdim A\otimes  U\}$$
where $U$ ranges over all division rings.
\item
The {\it third modification} of $\htr$ is defined to be
$$\thtr A=\sup\{\injdim A\otimes  U\}$$
where $U$ ranges over all division rings such that $A\otimes U$
is noetherian. 
\item
A simple artinian ring $S$ is called {\it homologically uniform} if 
$$\htr S =\fhtr S =\shtr S <\infty.$$
\item
A simple artinian ring $S$ is called {\it smooth} if 
$\pdim_{S^e} S<\infty$.
\end{enumerate}
If $A$ is a graded ring, then the graded version of
(a), (b), (c) and (d) can be defined and are denoted by $\htrgr_{\gr}$, 
$\fhtrgr_{\gr}$, $\shtrgr_{\gr}$ and $\thtrgr_{\gr}$ respectively.
\end{definition}

These definitions were implicitly suggested by the work of Resco 
\cite{Re1, Re2, Re3}, Rosenberg \cite{Ro}, Schofield \cite{Sc1} 
and Stafford \cite{St}. 
The idea in Definition \ref{xxdefn1.1}(e,f) of working in the class 
of simple artinian algebras instead of division algebras is due 
to Schofield \cite{Sc1}. Smooth simple artinian rings 
are called {\it regular} by Schofield \cite{Sc1} (see also
Lemma \ref{xxlem1.3}). 

We are mainly interested in $\htr$, but the modifications $\fhtr$,
$\shtr$ and $\thtr$ are closely related to $\htr$. In fact we are 
wondering whether division rings of finite $\htr$ are always 
homologically uniform [Question \ref{xxque7.9}] . 
For several classes of division algebras, $\fhtr$ and $\shtr$ are 
relatively easy to compute; and then $\htr$ is computable by 
the following easy lemma.

An algebra $A$ is called {\it doubly noetherian} if $A^\e$ is 
noetherian, and it is called {\it rationally noetherian} if $A\otimes U$ is 
noetherian for every division ring $U$.

\begin{lemma}
\label{xxlem1.2} 
Let $S$ be a simple artinian algebra. 
\begin{enumerate}
\item 
$\fhtr S \leq \htr S \leq \shtr S$. If $\fhtr S\geq \shtr S$ and 
$\shtr S<\infty$, 
then $S$ is homologically uniform. 
\item
$\thtr S\leq \shtr S$. If $S$ is rationally noetherian, then 
equality holds.
\item
If $S$ is doubly noetherian, then $\htr S\leq \thtr S$.
\end{enumerate}
\end{lemma}

By Proposition 
\ref{xxprop7.1} below there is a smooth, homologically uniform, 
but not doubly noetherian, commutative field $F$ over $k$ such that 
$\thtr F<\htr F$.

\begin{lemma}
\label{xxlem1.3} Let $S$ be a simple artinian ring.
In parts(a,c) suppose $S$ is smooth and $\pdim_{S^e} S
=n<\infty$. 
\begin{enumerate}
\item \cite[Lemma 2, p.\ 269]{Sc1}
$\gldim S^\e=n\geq \gldim S\otimes U$
for every simple artinian ring $U$. As a consequence, $\shtr S\leq n$.
\item
$S$ is smooth if and only if $\gldim S^e<\infty$. In
this case $\pdim_{S^e} S=\gldim S^e$.
\item
If $\Ext^n_{S^\e}(S,\boplus_I S^\e)\cong
\boplus_I \Ext^n_{S^\e}(S,S^\e)$ for any index set $I$ \textup{(}e.g.,
if $S$ is doubly noetherian\textup{)}, then $S$ is homologically uniform
and $\htr S=n$.
\end{enumerate}
\end{lemma}

\begin{proof} 
(b) Follows from (a).

(c) We claim that $\fhtr S=n$, which is equivalent to 
$\Ext^n_{S^\e}(S,S^\e)\neq 0$. Since $\pdim_{S^\e} S=n$,
there is an $S^\e$-module $M$ such that $\Ext^n_{S^\e}(S,M)
\neq 0$ and $\Ext^{n+1}_{S^\e}(S,-)=0$. There is an index set $I$
and a short exact sequence of $S^\e$-modules
$$0\to N\to \boplus_{I} S^\e\to M\to 0.$$
Applying $\Ext^i_{S^\e}(S,-)$ to the above exact sequence we have a long
exact sequence
$$\to \Ext^n_{S^\e}(S,N)\to \Ext^n_{S^\e}(S,\boplus_{I} S^\e)\to 
\Ext^n_{S^\e}(S,M)\to \Ext^{n+1}_{S^\e}(S,N)\to.$$
Since $\Ext^{n+1}_{S^\e}(S,N)=0$ and $\Ext^n_{S^\e}(S,M)\neq 0$, the 
above exact sequence implies that $\Ext^n_{S^\e}(S,\boplus_{I} S^\e)
\neq 0$. By hypothesis, $\Ext^n_{S^\e}(S,\boplus_I S^\e)\cong
\boplus_I \Ext^n_{S^\e}(S,S^\e)$, hence we have $\Ext^n_{S^\e}(S,S^\e)
\neq 0$. Thus our claim is proved. The assertion follows from part 
(a) and  Lemma \ref{xxlem1.2}(a).
\end{proof}

Lemma \ref{xxlem1.3} says that every doubly noetherian smooth 
division algebra $D$ is homologically uniform and $\htr D
=\gldim D^\e$. We can use Lemma \ref{xxlem1.3} and results 
of Resco and Stafford to compute the $\htr$ of some division 
rings. For example, the commutative field $D=k(x_1,\cdots,x_n)$ 
is homologically uniform with $\htr=n$, since $\gldim D^\e=n$ 
\cite[Theorem p.\ 215]{Re3}. A similar statement holds for the 
quotient division rings of $U(\mathfrak g)$ and $kG$ (see details 
in \cite[Theorem, p.\ 33]{St}). Proposition \ref{xxprop1.8}
below is also useful for such a  computation. 

Our main result Theorem \ref{xxthm0.2} deals with the case 
when $Q(A)$ may fail to be doubly noetherian, and Theorem \ref{xxthm0.3}
deals with the case when $Q(A)$ may fail to be smooth.

Let us now review the basic properties of the classical
transcendence degree of commutative fields over $k$. Let 
$F\subset G$ be commutative fields over $k$.

\begin{enumerate}
\item[(TD1)]
$\tr k(x_1,\cdots,x_n)=n$ for every $n\geq 0$.
\item[(TD2)]
$\tr F\leq \tr G$. 
\item[(TD3)]
If $\dim_F G$ is finite, then $\tr F =\tr G$.
\item[(TD4)]
If $G=F(x)$, then $\tr G=\tr F+1$.
\item[(TD5)]
If $\{F_i\}$ is a directed set of subfields of $G$ such that
$G=\bigcup F_i$, then $\tr G=\sup\{\tr F_i\}$.
\item[(TD6)]
If $G$ is finitely generated as a field and $\tr F =\tr G$,
then $\dim_F G$ is finite.
\end{enumerate}

\bigskip

We will try to prove some versions of (TD1-TD4) for $\htr$. 
However, Proposition \ref{xxprop7.1}(b) below shows that (TD5) 
fails for $\htr$, which is an unfortunate deficiency of 
homological transcendence degree. And we have not proven any 
generalization of (TD6). The following lemma is a 
collection of some well-known facts.

\begin{lemma} 
\label{xxlem1.4}
Let $A\subset B$ be $k$-algebras.
\begin{enumerate}
\item
Assume that $B=P\boplus N$ as left $A$-modules with $P$ a projective
generator of the category of left $A$-modules and that a similar 
decomposition holds for the right $A$-module $B$. Suppose the right
$A$-module $B$ is flat. Then $\injdim A\leq \injdim B$.
\item
If $A$ is an $A$-bimodule direct summand of $B$, then $\gldim A\leq 
\gldim B+\pdim_A B$.
\item
If $B=A[x;\alpha,\delta]$ where $\alpha$ is an automorphism
and $\delta$ is an $\alpha$-derivation, then 
$\injdim A\leq \injdim B\leq \injdim A+1$. 
\item
If $B$ is a localization of $A$, then $\gldim B\leq \gldim A$.
\item
If $A$ is noetherian and $B$ is a localization of $A$, then 
$\injdim B\leq \injdim A$.
\end{enumerate}
\end{lemma}

\begin{proof} (a) Since $B$ is a flat $A^\circ$-module, the 
$\Hom$-$\otimes$ adjunction gives the isomorphism
\begin{equation}
\label{E1.4.1}
\Ext^i_B(B\otimes_A M,B)\cong \Ext^i_A(M,B)
\tag{E1.4.1}
\end{equation}
for all $A$-modules $M$ and all $i$. The isomorphism \eqref{E1.4.1} is 
also given in Lemma \ref{xxlem2.2}(b) below. Since the right $A$-module 
$B$ contains a projective generator as a direct summand, $B\otimes_A M
\neq 0$ for every $M\neq 0$. So \eqref{E1.4.1} implies that 
$\injdim_B B\geq \injdim_A B$. Since the left $A$-module $B$ contains 
a projective generator as a direct summand, $\injdim_A B\geq 
\injdim_A A$. The assertion follows.

(b) This is \cite[Theorem 7.2.8]{MR}

(c) By part (a), $\injdim A\leq \injdim B$. 

By \cite[Proposition 7.5.2]{MR} (or \cite[Lemma 1, p.\ 268]{Sc1}),
for any $B$-module $M$, there is an exact sequence
$$0\to B\otimes_A (^\alpha M)\to B\otimes_A M\to M\to 0$$
of $B$-modules. This short exact sequence induces a long exact 
sequence
$$\cdots\to \Ext^i_B(B\otimes_A (^\alpha M),B)\to 
\Ext^{i+1}_B(M,B)\to \Ext^{i+1}_B(B\otimes_A M,B)
\to \cdots.$$
Since $B$ is a flat $A^\circ$-module, the $\Hom$-$\otimes$ adjunction 
gives 
$$\Ext^{i}_B(B\otimes_A M, B)\cong \Ext^i_A(M, B)=0$$
for all $i>\injdim A$ and all $M$. Thus the two ends of the above 
long exact sequence are zero, which implies that the middle term 
$\Ext^{i+1}_B(M,B)=0$ for all $i>\injdim A$. This shows that 
$\injdim B\leq \injdim A+1$.

(d) This is \cite[Corollary 7.4.3]{MR}.

(e) This is also well-known, and is a special case of Lemma 
\ref{xxlem2.3} below.
\end{proof}

\begin{proposition}
\label{xxprop1.5} Let $D\subset Q$ be simple artinian algebras.
\begin{enumerate}
\item
$\htr D \leq \htr Q$.
\item
If $A$ is Morita equivalent to $D$, then $\htr D =\htr A$.
\item
If $Q$ is finite over $D$ on the left, or on the right, 
then $\htr D =\htr Q $.
\item
If $D$ is PI and its center $C$ is finitely generated
over $k$ as a field, 
then $\htr D=\tr C=\tr D$.
\item
$\htr D =\htr D^\circ $.
\end{enumerate}
\end{proposition}

\begin{remark}
\label{xxrem1.6}
\begin{enumerate}
\item
The hypothesis ``$C$ is finitely generated
over $k$ as a field'' in part (d) of Proposition 
\ref{xxprop1.5} is necessary as Proposition
\ref{xxprop7.1} shows. 
\item
There are division algebras $D$ such that $D\not\cong D^\circ$.
For example, let $F$ be a field extension of $k$ such that
the Brauer group of $F$ has an element $[D]$ of order 
larger than $2$. Then the central division ring $D$ 
corresponding to $[D]$ has the property 
$D\not\cong D^\circ$.
If we want such a division algebra that is infinite over 
its center, then let $Q$ be the Goldie quotient ring
of $D\otimes k_q[x,y]$ where $q$ is not a root of 1. 
It is eay to check that the center of $Q$ is $F\otimes 
k\cong F$ and $Q$ is infinite over $F$. 
Since $D$ is the division subring of $Q$ 
consisting of all elements integral over the center $F$, 
then $Q\not\cong Q^{\circ}$.
\end{enumerate}
\end{remark}

\begin{proof}[Proof of Proposition \ref{xxprop1.5}]
(a) Since $D$ is simple artinian, every nonzero (left or right) 
$D$-module is a projective generator. So the $D$-module $Q$ and 
the $D^\circ$-module $Q^\circ$ are projective generators. Hence 
$Q\otimes Q^\circ$ is a projective generator over $D\otimes D^{\circ}$. 
Similarly, $Q\otimes Q^\circ$ is a  projective generator over 
$D\otimes D^{\circ}$ on the right. The assertion follows from
Lemma \ref{xxlem1.4}(a). 

(b) Since $A$ and $D$ are Morita equivalent and both are simple
artinian, $A\cong \mathrm{M}_n(B)$ and $D\cong \mathrm{M}_s(B)$ for some
division algebra $B$ and some $n,s$. So we may assume that $A$
is a division ring and $D=\mathrm{M}_n(A)$. Now $D^\e=
(\mathrm{M}_n(A))^\e\cong \mathrm{M}_{n^2}(A^\e)$. Hence we have 
$\injdim D^\e=\injdim \mathrm{M}_{n^2}(A^\e)=\injdim A^\e$.

(c) By part (a) it suffices to show that $\htr Q\leq \htr D$. Assume
the right $D$-module $Q$ is finite. Then $B:=\End_{D^\circ} (Q)$ is 
Morita equivalent to $D$. Also by the definition of $B$ there is a natural
injection $Q\to B$. By (a,b), $\htr Q\leq \htr B
=\htr D$. By symmetry the assertion holds when $Q$ is finite over
$D$ on the left. 

(d) Since $D$ is PI, $D$ is finite over its center $C$ 
\cite[Theorem 13.3.8]{MR}.
By part (c), $\htr D =\htr C$. It suffices to show that $\htr C
=\tr C$. Let $F$ be a subfield of $C$ such that 
$F\cong k(x_1,\cdots,x_n)$ for some integer $n$ and that 
$C$ is algebraic over $F$. Then $C$ is finite over $F$ and 
$\tr C=\tr F=n$. By part (c) it suffices to show that 
$\htr F =n$. By Lemma \ref{xxlem1.3} and \cite[Theorem, p.\ 215]{Re3}, 
$\htr F=n$. Hence the assertion follows.

(e) This follows from the fact that there is an anti-automorphism
$D^\e\to D^\e$.
\end{proof}

Similarly one can prove the following version of 
Proposition \ref{xxprop1.5} for $\shtr$. 

\begin{proposition}
\label{xxprop1.7} 
Let $D\subset Q$ be simple artinian algebras.
\begin{enumerate}
\item
$\shtr D\leq \shtr Q$.
\item
If $Q$ is finite over $D$ on the left (or on the right), 
then $\shtr D=\shtr Q$. As a consequence, if a simple artinian ring 
$A$ is Morita equivalent to $D$, then $\shtr D=\shtr A$.
\item
If $D$ is PI and its center, denoted by $C$, is finitely 
generated as a field, then $\shtr D=\tr C$.
\end{enumerate}
\end{proposition}

Unlike Proposition \ref{xxprop1.5}(e), we don't know whether
$\shtr D =\shtr D^\circ$ or not.

Recall from \cite{Sc1} that a simple artinian ring is {\it stratiform 
over $k$} if there is a chain of simple artinian rings
$$S=S_n\supset S_{n-1}\supset \cdots \supset S_1\supset S_0=k$$
where, for every $i$, either (i) $S_{i+1}$ is finite over $S_i$ on both
sides; or (ii) $S_{i+1}$ is isomorphic to $S_i(x_i;\alpha_i,\delta_i)$
for an automorphism $\alpha_i$ of $S_i$ and and $\alpha_i$-derivation
$\delta_i$ of $S_i$. Such a chain of simple artinian rings is called a 
{\it stratification} of $S$. The {\it stratiform length} of $S$ is
the number of steps in the chain that are of type (ii). 
One basic property proved in \cite{Sc1} is that the stratiform 
length is an invariant of $S$.

\begin{proposition}
\label{xxprop1.8} If $S$ is a stratiform simple artinian ring
of stratiform length $m$, then $S$ is rationally noetherian, homologically 
uniform and $\htr S =m$.
\end{proposition}

\begin{proof} It follows from induction on the steps of the stratification 
that $S$ and $S^\circ$ are rationally noetherian (and hence doubly 
noetherian). But $S$ might not be smooth. 

Next we show $\fhtr S=m$. Applying \cite[Lemma 20, p.\ 277]{Sc1} 
to the $S$-bimodule $S$, we have the following statement:
there is a simple artinian ring $S''\supset S$ such that 
$S''$, as $(S''\otimes S^\circ)$-module, has projective
dimension $m$. Since $S^{\circ}$ is rationally noetherian, 
$S''\otimes S^\circ$ is noetherian. Hence we have  
$$\Ext^i_{S''\otimes S^\circ}(S'',S''\otimes S^\circ)=
\begin{cases} {\text {nonzero}} & {\text {if}} \quad i=m,\\
0 & {\text {if}} \quad i>m.
\end{cases}$$
Since $S''\otimes S^\circ$ is projective over $S^\e$, we have 
$$\Ext^i_{S''\otimes S^\circ}(S'',S''\otimes S^\circ)
\cong \Ext^i_{S^\e}(S,S''\otimes S^\circ)$$ which is a direct sum of 
copies of $\Ext^i_{S^\e}(S,S^\e)$. Thus 
$$\Ext^i_{S^\e}(S,S^\e)=
\begin{cases} {\text {nonzero}} & {\text {if}} \quad i=m,\\
0 & {\text {if}} \quad i>m.
\end{cases}$$
Therefore $\fhtr S=m$. 

By Lemma \ref{xxlem1.2}(a) it remains to show that $\shtr S\leq m$.
We use induction on the steps of the stratification. Suppose 
$\shtr S_{n-1}$ is no more than the stratiform length of 
$S_{n-1}$. We want to show that this statement holds for $S_n$.

Case (i): $S_n$ is finite over $S_{n-1}$ on both sides. By 
Proposition \ref{xxprop1.7}(b), $\shtr S_{n}=\shtr S_{n-1}$.
The claim follows.

Case (ii): $S_n=S_{n-1}(x;\alpha,\delta)$. Let $U$ be any simple 
artinian ring and let $A=S_{n-1}\otimes U$ and 
$B=S_{n-1}[x;\alpha,\delta]\otimes U=A[x;\alpha,\delta]$.
By Lemma \ref{xxlem1.4}(c), $\injdim B\leq \injdim A+1$.
Since $S_n\otimes U$ is a localization of the noetherian 
ring $B$, we have $\injdim S_n\otimes U\leq \injdim B$ by
Lemma \ref{xxlem1.4}(e). Combining these two inequalities, the
claim follows.
\end{proof}

Next we give a list of known examples, and a few more examples will 
be given in Section \ref{xxsec7}. 

\begin{example}
\label{xxex1.9} 
\begin{enumerate}
\item Let $F$ be a separable field extension of $k$ that is
finitely generated as a field. Then $F$ is rationally noetherian and 
$\gldim F^\e=\tr F<\infty$. Hence $F$ is smooth, homologically 
uniform and $\htr F=\tr F$.
\item
Let $F$ be the commutative field $k(x_1,x_2,\ldots)$, which is an 
infinite pure transcendental extension of $k$. The ring $F$ is not 
doubly noetherian. For each integer $m$ let $F_m$ be the subfield 
$k(x_1,\ldots,x_m)\subset F$. Then $F_m$ is rationally noetherian,
smooth, homologically uniform with $\htr$ $m$. Since $F_m\subset F$
for all $m$, one sees that 
$$\htr F=\shtr F=\thtr F=\infty.$$
But $\fhtr F=-\infty$ since $\Ext^i_{F^\e}(F,F^\e)=0$ for all
$i$ \cite[Example 3.13]{YZ1}. 
\item
Let $D$ be a simple artinian ring finite dimensional over $k$. By 
Proposition \ref{xxprop1.5}(c), $\htr D=\htr k=0$.
\item
Let $F$ be a finite dimensional purely inseparable field extension of $k$.
Since $\gldim F^\e=\infty$, $F$ is not smooth over $k$.
By part (c) $\injdim F^\e=\htr F=0$.
\item
Let $D_n$ be the $n$-th Weyl skew field. Since $D_n$ is rationally 
noetherian and $\gldim D_n^\e=2n$ \cite{Ro,St}, by Lemma \ref{xxlem1.3}, 
$D_n$ is smooth and homologically uniform and $\htr D_n=2n$. 
\item
Let $D(\mathfrak g)$ be the quotient division ring of 
the universal enveloping algebra $U(\mathfrak g)$ of 
a finite dimensional Lie algebra $\mathfrak g$. Then
$D(\mathfrak g)$ is rationally noetherian and 
$\gldim D(\mathfrak g)^\e=\dim_k {\mathfrak g}$ \cite{St}. 
By Lemma \ref{xxlem1.3},  $D(\mathfrak g)$ is smooth and 
homologically uniform and $\htr D=\dim_k {\mathfrak g}$.
A similar statement holds for quotient division rings of
group rings $kG$ studied in \cite{St}.
\item
Let $\{p_{ij}\;|\; i<j\}$ be a set of nonzero scalars in $k$.
Let $A$ be the skew polynomial ring $k_{\{p_{ij}\}}[x_1,\ldots,x_n]$
that is generated by elements $x_1,\ldots,x_n$ and subject to the 
relations $x_jx_i=p_{ij}x_ix_j$ for all $i<j$. Let $Q$ be the
quotient division ring of $A$. Then $Q$ is a stratiform division
ring of stratiform length $n$. Hence $Q$ is rationally noetherian, 
homologically uniform, and $\htr Q=n$. This is a generalization of (TD1).
 Since $Q^\e$ is a localization of another skew polynomial ring of 
finite global dimension, $Q$ is smooth. 
\end{enumerate}
\end{example}

\section{Polynomial extension}
\label{xxsec2}

In this section we discuss the property (TD4) for $\htr$. We have 
not yet proved a satisfactory generalization of  (TD4). 
Let $S$ be a simple artinian ring with automorphism $\alpha$
and $\alpha$-derivation $\delta$ of $S$. The Goldie quotient 
ring of $S[t;\alpha,\delta]$ is denoted by $S(t;\alpha,\delta)$. 
We don't know if 
$$\htr S(t;\alpha,\delta)=\htr S+1$$
holds in general, but we present some partial results in Proposition
\ref{xxprop2.7} below.

Recall that the third modification of $\htr$ is
$$\thtr S=\sup\{\injdim S\otimes U\}$$
where $U$ ranges over all division algebras such that 
$S\otimes U$ is noetherian. If $S$ is a doubly noetherian simple artinian 
ring, then 
$$\fhtr S\leq \htr S \leq \thtr S\leq \shtr S.$$
For doubly noetherian simple artinian rings $S$, $\thtr$ is a good
replacement for $\shtr$. In this case we call $S$ {\it weakly uniform} 
if 
$$\fhtr S=\htr S=\thtr S<\infty.$$ 
If $S$ is rationally noetherian, then ``weakly
uniform'' is equivalent to ``homologically uniform''. 

\begin{lemma}
\label{xxlem2.1} Let $S$ be a simple artinian ring and
let $Q=S(t;\alpha,\delta)$.
\begin{enumerate}
\item
$S$ is doubly noetherian if and only if $Q$ is.
\item
$S$ is rationally noetherian if and only if $Q$ is.
\item
If $S$ is smooth, so is $Q$. The converse holds when
$\delta=0$. 
\end{enumerate}
\end{lemma}

\begin{proof} Note that $Q^\circ\cong 
S^\circ(t;\alpha^{-1},-\delta\alpha^{-1})$. 

(a) If $S^\e=S\otimes S^{\circ}$ is noetherian, so is 
$S[t;\alpha,\delta]\otimes S^\circ[t;\alpha^{-1},
-\delta\alpha^{-1}]$. Therefore its localization $Q\otimes Q^\circ$ 
is noetherian.

In the other direction, we suppose $Q^\e$ is noetherian.
Since $Q$ is faithfully flat (and projective) as left and
right $S$-module, $Q^\e$ is a faithfully flat left module over 
$S^\e$. Hence $S^\e$ is left (and hence right) noetherian.

(b) Similar to part (a).

(c) If $S$ is smooth, an argument similar to the proof of part (a)
shows that $Q$ is smooth. 

To show the converse we assume that $\delta=0$. Decompose $Q$ into 
$Q=S\boplus C$ where 
$$C=\{f(t)(g(t))^{-1} \;|\; \deg_t f(t)<\deg_t g(t)\}\oplus
(\boplus_{n=1}^{\infty}t^n S).$$ 
Hence $S$ is a $S$-bimodule 
direct summand of $Q$. Thus $S^\e$ is a $S^\e$-bimodule direct
summand of $Q^\e$. The assertion follows from Lemma \ref{xxlem1.4}(b).
\end{proof}

It is not clear to us if the ``converse'' part of Lemma 
\ref{xxlem2.1}(c) holds when $\delta\neq 0$.

The following lemma is basically \cite[Lemma 3.7]{YZ2}. Note that in 
\cite[Lemma 3.7]{YZ2}, an extra hypothesis ``$M$ being bounded below'' 
was forgotten. Various versions of the following lemma exist in 
the literature, especially for modules instead of complexes. 

Let $\Mod A$ denote the category of $A$-modules and let $\D(\Mod A)$
denote the derived category of $\Mod A$. If $A$ is graded, 
$\GrMod A$ is the category of graded $A$-modules and  
$\D(\GrMod A)$ is the derived category of $\GrMod A$. We refer 
to \cite{Ye1} for basic material about complexes and derived categories.

A complex $L\in \D^{-}(\Mod A)$ is called {\it pseudo-coherent} 
if $L$ has a bounded above free resolution 
$P=(\cdots \to P^i\to P^{i+1}\to \cdots)$ such that each component 
$P^i$ is a finite free $A$-module 
\cite[Expos\'e 1 (L. Illusie), Section 2]{SGA6}. 
If $L$ is pseudo-coherent, then 
$$\RHom_A(L,\boplus_{i\in I} M_i)\cong \boplus_{i\in I}\RHom_A(L,M_i)$$
and
$$\Ext^n_A(L,\boplus_{i\in I} M_i)\cong \boplus_{i\in I}\Ext^n_A(L,M_i)
\quad
\text{for all $n$}$$ 
where $\{M_i\}_{i\in I}$ is a set of uniformly bounded below complexes. 

There are two different definitions of injective dimension 
of a complex existing in the literature, one of which is given 
as follows. Let $X$ be a bounded below complex of $A$-modules. 
Then the {\it injective dimension} of $X$ is defined to be
$$\injdim_A X=\sup\{i\;|\; Y^i\neq 0\}$$
where $Y$ is a minimal injective resolution of $X$. If $A$ is
$\mathbb Z$-graded, the graded $\injdim$ can be 
defined. If $\injdim_A X=n$, then $\Ext^i_A(M,X)=0$ for all $A$-modules
$M$ and for all $i>n$; and there is an $A$-module
$M$ such that $\Ext^n_A(M,X)\neq 0$. 

\begin{lemma} \cite[Lemma 3.7]{YZ2}
\label{xxlem2.2}
Let $A$, $B$ be algebras. Let $L$ be a complex in $\D^{-}(\Mod A)$.
\begin{enumerate}
\item Let $N$ be a $B$-module of finite flat dimension 
and let $M \in \D^{+}(\Mod A \otimes B^{\circ})$. 
Suppose $L$ is pseudo-coherent. Then the 
functorial morphism
$$\RHom_{A}(L, M) \otimes^{\rm L}_{B} N \to 
\RHom_{A}(L, M \otimes^{\rm L}_{B} N)$$
is an isomorphism in $\D(\Mod k)$.
\item Suppose $A \to B$ is a ring homomorphism such that 
$B$ is a flat $A^{\circ}$-module. Let $M \in \D^{+}(\Mod B)$.
Then the functorial morphism
$$\RHom_{A}(L, M) \to \RHom_{B}(B\otimes_A L, M)$$
is an isomorphism in $\D(\Mod k)$. 
\end{enumerate}
\end{lemma}

The following lemma is well-known; and  follows easily from
the above lemma.

\begin{lemma}
\label{xxlem2.3}
Let $A$ be a noetherian ring and let $B$ be any ring. Suppose $R$ is 
a bounded complex of $A\otimes B^\circ$-modules. Suppose that $A'$ and $B'$
are Ore localizations of $A$ and $B$ respectively and that $A'\otimes_A R
\cong R\otimes_B B'$ in $\D(\Mod A\otimes B^\circ)$. Then $\injdim_{A'} 
(A'\otimes_A R) \leq \injdim_A R$.
\end{lemma}

The following lemma is similar to \cite[Theorem 8, p.\ 272]{Sc1} and 
is known to many researchers. Note that there is a typographical error 
in the statement of \cite[Theorem 8, p.\ 272]{Sc1}: ``$i\neq 1$'' 
should be ``$i\geq 1$''.

\begin{lemma} 
\label{xxlem2.4}
Let $A$ be a left noetherian ring with an automorphism $\alpha$. 
Let $T=A[t^{\pm 1};\alpha]$. If $M$ is a $T$-module that is 
finitely generated as $A$-module, then 
$$\Ext^i_T(M,T)\cong \Ext^{i-1}_A({^\alpha M},A)$$
as $A^\circ$-modules, for all $i\geq 1$. 
\end{lemma}

\begin{definition}
\label{xxdefn2.5} A simple artinian algebra $S$ is called
{\it rigid} if $\RHom_{S^\e}(S,S^\e)\cong S^\sigma[-n]$ 
for some integer $n$ and some automorphism $\sigma$ of $S$;
or equivalently
$$\Ext^i_{S^\e}(S,S^\e)=\begin{cases} S^\sigma & 
{\text {if}} \quad i=n,\\
0 & {\text {if}} \quad i\neq n.
\end{cases}$$
\end{definition}

In the above definition $[n]$ denotes the $n$-th complex
shift and the bimodule $S^\sigma$ is defined by
$$a*s*b=as\sigma(b)$$
for all $a,s,b\in S$. Clearly $\fhtr S=n$. 

For example the $n$-th Weyl skew field $D_n$ in Example 
\ref{xxex1.9}(e) is rigid. This follows from the computation 
given at the end of \cite[Section 6]{YZ3}:
$$\Ext^i_{D_n^\e}(D_n,D_n^\e)=\begin{cases} D_n & 
{\text {if}} \quad i=2n,\\
0 & {\text {if}} \quad i\neq 2n.
\end{cases}$$
Other doubly noetherian division rings in Example \ref{xxex1.9}
are also rigid [Corollary \ref{xxcor6.11}]. If $S$ is not
doubly noetherian, then $S$ might not be rigid even if it is
homologically uniform [Proposition \ref{xxprop7.1}(d)].

\begin{proposition}
\label{xxprop2.6} Let $D$ be a doubly noetherian simple artinian 
ring and $B$ be the skew Laurent polynomial ring $D[t^{\pm 1};\alpha]$ 
where $\alpha$ is an automorphism of $D$. Let $Q=D(t;\alpha)$.
Then $\fhtr B=\fhtr Q=\fhtr D+1$. Furthermore $D$ is rigid if and only 
if $Q$ is rigid. 
\end{proposition}

\begin{proof} It is easy to reduce to the case when $D$ is a 
division ring, so we assume that $D$ is a division ring in the 
proof below. Since $D^\e$ is noetherian, so are $B^\e$ and $Q^\e$. 

We can view $B$ as a 
$\mathbb Z$-graded ring with $\deg t=1$ and $\deg D=0$.
Hence $B^\e$ is also $\mathbb Z$-graded and 
$B^\e\cong (D^\e[w^{\pm 1};\alpha\otimes\alpha])
[(t\otimes 1)^{\pm 1},\sigma]$ where $w=t\otimes t^{-1}$ and 
$\sigma: w\mapsto w, d_1\otimes d_2\mapsto \alpha(d_1)\otimes d_2$
for all $d_1\in D$ and $d_2\in D^\circ$. 
Clearly $B^\e$ is strongly $\mathbb Z$-graded.

Since $B^\e$ is noetherian and $B$ is a finite graded $B^\e$-module, 
$\Ext^i_{B^\e}(B,B^\e)$ can be computed in the category of 
the $\mathbb Z$-graded $B^\e$-modules. Since $B^\e$ is strong
$\mathbb Z$-graded, $\GrMod B^\e\cong \Mod C$ where 
$C=D^\e[w^{\pm 1};\alpha\otimes\alpha]$. Since the degree zero 
parts of $B^e$ and $B$ are equal to $C$ and $D$ respectively, we have
$\Ext^i_B(B,B^\e)\neq 0$ if and only if $\Ext^i_C(D,C)\neq 0$
where $D$ is a left $C$-module. Furthermore the degree zero part of
$\Ext^i_B(B,B^\e)$ is isomorphic to $\Ext^i_C(D,C)$. By Lemma 
\ref{xxlem2.4}, for all $i\geq 1$, 
\begin{equation}
\label{E2.6.1}
\Ext^{i}_C(D,C)\cong \Ext^{i-1}_{D^\e}(^{(\alpha\otimes\alpha)}D, 
D^\e)\cong \Ext^{i-1}_{D^\e}(D,D^\e)^{\alpha^{-1}\otimes \alpha^{-1}}.
\tag{E2.6.1}
\end{equation}
This shows that $\fhtr B= \fhtr D+1$.

Let $n$ be any integer such that $V:=\Ext^n_{B^\e}(B,B^\e)\neq 0$;
and let $W=\Ext^{n}_C(D,C)$. There is a natural $D$-bimodule structure
on $W$. Let's think about the left $D$-action on $W$. Since $D$ 
is a division ring, $W$ is a faithful $D$-module. Since 
$V$ is basically equal to $W[t^{\pm 1}]$, it is a faithful 
$B(=D[t^{\pm 1};\alpha])$-module. Therefore $Q\otimes_B V\neq 0$.
Similarly $V$ is a faithful $B^\circ$-module. Hence 
$Q\otimes_B V$ is a faithful $B^\circ$-module and 
$Q\otimes_B V\otimes_B Q\neq 0$. Since $B^\e$ is noetherian,
by Lemma \ref{xxlem2.2}, 
$$\Ext^n_{Q^\e}(Q,Q^\e)\cong \Ext^n_{B^\e}(B,Q^\e)\cong
V\otimes_{B^\e} Q^\e\cong Q\otimes_B V\otimes_B Q\neq 0.$$
This implies that $\fhtr B=\fhtr Q$.

If $\Ext^{n-1}_{D^\e}(D,D^\e)\cong D^\sigma$ for some automorphism
$\sigma$ of $D$, then by \eqref{E2.6.1} we have  $\Ext^{n}_C(D,C)
\cong D^{\sigma'}$ for another automorphism $\sigma'$ of $D$. The 
above argument shows that $\Ext^n_{Q^\e}(Q,Q^\e)\cong Q^{\sigma''}$ 
for some automorphism $\sigma''$ of $Q$. Therefore if $D$ is rigid 
so it $Q$. The converse can be proved similarly.  
\end{proof}

\begin{proposition}
\label{xxprop2.7}
Suppose $S$ is a doubly noetherian simple artinian ring. Let $Q=S(t;\alpha)$.
\begin{enumerate}
\item
If $S$ (or $Q$) is smooth, then $\htr Q=\htr S+1$.
\item
If $S$ is weakly uniform, then so is $Q$, and 
$\htr Q=\htr S+1$.
\item
If $\alpha=id_S$ and if $\fhtr Q=\htr Q$,  
then $\fhtr S=\htr S=\htr Q-1$. 
\end{enumerate}
\end{proposition}

\begin{proof} (a) The assertion follows from Lemma 
\ref{xxlem1.3} and Proposition \ref{xxprop2.6}.

(b) Note that $S\otimes U$ is noetherian if and only if
$Q\otimes U$ is. By Lemma \ref{xxlem1.4}(c,e), we have
$\thtr Q \leq \thtr S+1$. If $S$ is weakly uniform, we have
$$\thtr S+1=\fhtr S+1=\fhtr Q$$
where the last equality is Proposition \ref{xxprop2.6}. 
Combining these facts with Lemma \ref{xxlem1.2}(a,c),
we obtain that $Q$ is weakly uniform and that 
$\htr Q=\htr S+1$.

(c) By hypothesis and Proposition \ref{xxprop2.6} we have
$$\htr Q=\fhtr Q=\fhtr S+1\leq \htr S+1.$$
Let $A=S(t)\otimes_{k(t)}S^\circ(t)$. Then 
$A\cong Q^\e/(t\otimes 1-1\otimes t)$
where $(t\otimes 1-1\otimes t)$ is a central regular element of 
$Q^\e$. By Rees' lemma, $\injdim Q^\e\geq \injdim A +1$. We can decompose
$A$ into a direct sum of $S^\e$-modules
$$A=S^\e\boplus (\boplus_{n\geq 1}S^\e t^n)\boplus 
\{f_1g_1^{-1}\otimes f_2g_2^{-1}\;|\;\deg_t f_i<\deg_t g_i\}$$
$$\boplus \{S\otimes f_2g_2^{-1}\;|\;\deg_t f_2<\deg_t g_2\}
\boplus \{f_1g_1^{-1}\otimes S^\circ\;|\;\deg_t f_1<\deg_t g_1\}.
$$
Also, $A$ is flat over $S^\e$ on the right. By Lemma \ref{xxlem1.4}(a), 
$\injdim A\geq \injdim S^\e$. Hence we have
$$\htr Q=\injdim Q^\e\geq \injdim A+1\geq \injdim S^\e+1=\htr S+1.$$
Combining these inequalities with $\fhtr Q=\htr Q$, $\htr S+1\geq 
\fhtr S+1$ (follows from Lemma \ref{xxlem1.2}(a)) and 
$\fhtr Q=\fhtr S+1$ (follows from Proposition 
\ref{xxprop2.6}), one sees that all inequalities are equalities;
and hence $\htr S+1=\fhtr S+1=\htr Q$. 
\end{proof}

Now we are ready to prove Proposition \ref{xxprop0.4}.

\begin{proof}[Proof of Proposition \ref{xxprop0.4}]
Parts (a,b,d,e) are proved in Proposition \ref{xxprop1.5}
and part (c) in Proposition \ref{xxprop2.7}(a).
\end{proof}

In the rest of the paper we will compute $\htr$ for various 
classes of division algebras that are not in Example \ref{xxex1.9}.

\section{Review of Dualizing complexes}
\label{xxsec3}

The dualizing complex is one of the main tools in the 
computation of homological transcendence degree. In this section 
we review several basic definitions related to dualizing complexes. 
We refer to \cite{VdB, Ye1, YZ1} for other details.  Some material 
about local duality will be reviewed in Section \ref{xxsec6}. 

\begin{definition} 
\label{xxdefn3.1} Let $A$ be an algebra. A 
complex $R\in \D^{\rm b}(\Mod A^\e)$ is called a {\it dualizing 
complex} over $A$ if it satisfies the following conditions:
\begin{enumerate}
\item
$R$ has finite injective dimension over $A$ and over $A^{\circ}$
respectively.
\item
$R$ is pseudo-coherent over $A$ and over $A^{\circ}$ respectively.
\item 
The canonical morphisms $A\to \RHom_A(R,R)$ and 
$A\to \RHom_{A^{\circ}}(R,R)$ are isomorphisms in $\D (\Mod A^\e)$.
\end{enumerate}
If $A$ is $\mathbb Z$-graded, a graded dualizing complex is defined similarly.
\end{definition}

If $A$ is noetherian (or more generally, coherent) then the definition 
agrees with \cite[Definition 3.3]{Ye1} (or \cite[Definition 1.1]{YZ1} 
for $A=B$). 

Let $R$ be a dualizing complex over a noetherian ring $A$ and let $M$ 
be a finite $A$-module. The {\it grade} of $M$ with respect to $R$ 
is defined to be
$$j_R(M) = \inf\, \{q \mid \Ext^q_A(M, R) \neq 0 \} .$$
The grade of an $A^{\circ}$-module is defined similarly.

\begin{definition} 
\label{xxdefn3.2}
\cite{Ye2, YZ1}
A dualizing complex $R$ over a noetherian ring $A$ is called 
{\it Auslander} if
\begin{enumerate}
\item 
For every finite $A$-module $M$, every integer $q$ and every 
finite $A^{\circ}$-submodule $N \subset \Ext^q_A(M, R)$ one has 
$j_R(N) \geq q$.
\item 
The same holds after exchanging $A$ and 
$A^{\circ}$. 
\end{enumerate}
\end{definition}

The {\it canonical dimension} of a finite $A$-module $M$ with respect 
to $R$ is defined to be
$$\Cdim  M = - j_R(M).$$

Let $R$ be a complex of $A^\e$-modules, viewed as a complex of $A$-bimodules. 
Let $R^{\circ}$ denote the ``opposite complex'' of $R$ which is 
defined as follows: as a complex of $k$-modules $R=R^{\circ}$ and 
the left and right $A^{\circ}$-module actions on $R^{\circ}$ is 
given by
$$a \ast r \ast b=bra$$
for all $a, b\in A^{\circ}$ and $r\in R^{\circ}(=R)$. If $R\in \D(\Mod A^\e)$
then $R^{\circ}\in \D(\Mod (A^\circ)^\e)$. Since $(A^{\circ})^\e$
is isomorphic to $A^\e$, there is a natural isomorphism $\D(\Mod A^\e)\cong
\D(\Mod (A^{\circ})^\e)$. The following definition is due 
to Van den Bergh \cite[Definition 8.1]{VdB}.

\begin{definition} \cite{VdB} 
\label{xxdefn3.3}
A dualizing complex $R$ over $A$ is called {\it rigid} if 
there is an isomorphism
$$\rho: R \to \RHom_{A^\e}(A, R \otimes R^{\circ})$$ 
in $\D(\Mod A^\e)$. Here the left $A^\e$-module structure of
$R\otimes R^{\circ}$ comes from the left $A$-module structure
of $R$ and the left $A^{\circ}$-module structure of $R^{\circ}$.
To be precise $(R, \rho)$ is called 
a {\it rigid} dualizing complex and the isomorphism $\rho$ is 
called a {\it rigidifying isomorphism}.
\end{definition}

A simple artinian ring $S$ is rigid (see Definition \ref{xxdefn2.5}) 
if and only if $S$ has a rigid dualizing complex. In fact 
an easy computation shows that 
$\RHom_{S^\e}(S,S^\e)\cong S^{\sigma}[-n]$ if and only if 
$R:=S^{\sigma^{-1}}[n]$ is a rigid dualizing complex over $S$. 

When $A$ is connected graded, there is a notion of balanced
dualizing complex introduced in \cite{Ye1}, which is 
related to the rigid dualizing complex. Let $A$ be a 
connected graded algebra and let $\fm=A_{>0}$. 
Let $\Gamma_{\fm}$ denote the $\fm$-torsion functor $\lim_{n\to \infty} 
\operatorname{Hom}_A(A/\fm^n, -)$ (also see Section \ref{xxsec6}). 
If $M$ is a graded
$A$-module, let $M'$ denote the graded vector space dual of $M$.

\begin{definition}\cite{Ye1} 
\label{xxdefn3.4} 
A graded dualizing complex $R\in \D^{\rm b}(\GrMod A^\e)$ over a connected 
graded ring $A$ is called {\it balanced} if there are isomorphisms 
$$\operatorname{R\Gamma}_{\fm}(R)\cong A'
\cong \operatorname{R\Gamma}_{\fm^\circ}(R)$$
in $\D^{b}(\GrMod A^\e)$.
\end{definition}

By \cite[Proposition 8.2(2)]{VdB} a balanced dualizing complex
over a noetherian connected graded ring is rigid after 
forgetting the grading.

\begin{definition}
\label{xxdefn3.5}
A connected graded ring $A$ is called {\it Artin-Schelter Gorenstein}
(or AS Gorenstein) if 
\begin{enumerate}
\item
$A$ has graded injective dimension $n<\infty$ on the left and on the right,
\item 
$\Ext^i_A(k,A)=\Ext^i_{A^\circ}(k,A)=0$ for all $i\neq n$, and 
\item
$\Ext^n_A(k,A)\cong \Ext^n_{A^\circ}(k,A)\cong k(l)$ for some $l$.
\end{enumerate}
If moreover $A$ has finite graded global dimension, then 
$A$ is called {\it Artin-Schelter regular} (or AS regular).
\end{definition}

In the above definition $(l)$ denotes the $l$th degree shift
of a graded module. If $A$ is AS regular, then $\gldim A=n=\injdim A$.
By \cite{Ye1} if $A$ is noetherian and AS Gorenstein (or regular),
then $A$ has a balanced dualizing complex $A^\sigma(-l)[-n]$ for some 
automorphism $\sigma$. Note that in Definition \ref{xxdefn3.5}
neither $A$ is noetherian nor is the GK-dimension of $A$ 
finite.

\section{Computation of $\fhtr$}
\label{xxsec4}

In this section we use Van den Bergh's rigidity formula to 
compute $\fhtr$ of some division algebras. 

Let $A$ be an algebra and let $S$ be a left and right Ore set of 
regular elements of $A$. Let $B=S^{-1}A=AS^{-1}$. An $A$-bimodule 
complex $R$ is called {\it evenly localizable} to $B$ if 
$$B\otimes_A R\to 
B\otimes_A R\otimes_A B\quad\text{ and }\quad
R\otimes_A B\to B\otimes_A R
\otimes_A B$$ are quasi-isomorphisms \cite[Definition 5.8]{YZ3}. 
If $B$ is $Q(A)$, the total Goldie quotient ring of $A$, then we 
simply say $R$ is evenly localizable without reference to $B$. 
It is easy to see that $R$ is evenly localizable to $B$ if and 
only if $\HB^i(R)$ is evenly localizable to $B$ for all $i$. 
The following lemma was proved a few times in slightly 
different versions (e.g., \cite[Theorem 6.2]{YZ3}).

\begin{lemma}
\label{xxlem4.1} Let $A$ be an algebra and let $B$ be a localization
of $A$ with respect to an Ore set $S$. Let $R$ be a dualizing 
complex over $A$. Assume that 
\begin{enumerate}
\item[(i)] 
$R$ is evenly localizable to $B$, and 
\item[(ii)]
either $A$ is noetherian or $B$ has finite global
dimension.
\end{enumerate}
Then $R_B:=B\otimes_A R\otimes_A B$ is a dualizing complex over $B$.

The graded version of the assertion also holds.
\end{lemma}

\begin{proof} We only sketch a proof in the case when $B$ 
has finite (ungraded) global dimension. By the definition of even 
localizibility we have $R\otimes_A B\cong R_B\cong B \otimes_A R$. 

To prove $R_B$ is a dualizing complex over $B$ we need to show 
(a,b,c) in Definition \ref{xxdefn3.1}. 
Part (a) is clear since $B$ has finite global dimension and 
$R_B$ is bounded. Part (b) follows from the fact
that the pseudo-coherence is preserved under flat change of rings. 
Part (c) follows from Lemma \ref{xxlem2.2} 
and the fact $R$ is pseudo-coherent.

The noetherian case is similar; in fact it was proved in 
\cite[Theorem 6.2(a)]{YZ3}. 
\end{proof}

Let $A$ be a Goldie prime ring and let $Q(A)$ denote the Goldie
quotient ring of $A$. Then $Q(A)$ is simple artinian; and in 
particular, it has global dimension 0. If $A$ is graded Goldie 
prime, let $Q_{\rm gr}(A)$ denote the graded Goldie quotient ring of 
$A$. Then $Q_{\rm gr}(A)$ is graded simple artinian of graded global 
dimension 0. As an ungraded ring, $Q_{\rm gr}(A)$ is noetherian and 
has global dimension at most 1. 

Suppose $R$ is a dualizing complex over $A$ that is evenly 
localizable to $Q:=Q(A)$. Since $Q$ is simple artinian of
global dimension 0, Lemma \ref{xxlem4.1} applies and $R_Q(\cong
R\otimes_A Q\cong Q\otimes_A R)$ is a dualizing complex 
over $Q$. Let $B$ be any simple artinian ring. Then
a dualizing complex over $B$ is isomorphic to $P[n]$ where 
$P$ is an invertible $B$-bimodule \cite[Theorem 0.2]{YZ5}; and 
every invertible $B$-bimodule is isomorphic to ${B^\sigma}$ for
some automorphism $\sigma$ of $B$. Hence every dualizing complex 
over $Q$ is isomorphic to ${Q^\sigma}[n]$ for some $n$ and some 
automorphism $\sigma$ of $Q$. Therefore one has
$$\{i \;|\; Q\otimes_A \HB^i(R)\neq 0\}=
\{i \;|\; \HB^i(R)\otimes_A Q\neq 0\}=
\{i \;|\; \HB^i(R_Q)\neq 0\}=\{-n\}.$$
Use this equation we define {\it hammerhead} of $R$ to be $n$
and write $\xi(R)=n$. In the graded setting, every graded 
dualizing complex over a graded simple artinian ring $Q_{\rm gr}$ 
is of the form ${Q_{\rm gr}^\sigma}(l)[n]$. So one can define a 
graded version of this nation, called the {\it hammerhead} of 
the graded dualizing complex $R$, and denoted by $\xi_{\gr}(R)$.

\begin{proposition}
\label{xxprop4.2} Let $A$ be a Goldie prime ring and let 
$Q=Q(A)$. Let $R$ be a rigid dualizing complex over $A$ that 
is evenly localizable to $Q$. If the $A^\e$-module $A$ is 
pseudo-coherent, then $Q$ is rigid and $\fhtr Q=\xi(R)$. 

If $A$ is graded, then the graded version of the assertion also 
holds. 

In the graded case, we have further $\fhtr Q=\fhtr_{\gr} Q_{\rm gr}(A)=
\xi(R)=\xi_{\gr}(R)$.
\end{proposition}

\begin{proof} It is easy to show that the ring $Q^\e$ is 
an Ore localization of $A^\e$. So $Q^\e$ is a flat $A^\e$-module.
By Lemma \ref{xxlem2.2}(a) and pseudo-coherence of $A$,
$$\RHom_{A^\e}(A,R\otimes R^{\circ})\otimes_{A^\e} Q^\e
\cong \RHom_{A^\e}(A,(R\otimes R^{\circ})\otimes_{A^\e} Q^\e)=:(*).$$
Since $R_Q:=Q\otimes_A R\otimes_A Q$ is a dualizing complex over $Q$, 
$R\otimes_A Q\cong X[n]$ and 
$R^{\circ}\otimes_A Q^{\circ}= X^{\circ}[n]$ where $X$ is
is isomorphic to a $Q$-bimodule $Q^\tau$ for some
automorphism $\tau$ and $n=\xi(R)$. Hence
$$(*)=\RHom_{A^\e}(A,X\otimes X^{\circ}[2n])
\cong \RHom_{Q^\e}(Q^\e\otimes_{A^\e}A,X\otimes X^{\circ}[2n])$$
where the last isomorphism is Lemma \ref{xxlem2.2}(b). Note that
$Q^\e\otimes_{A^\e}A\cong Q$ as $Q^\e$-module. By the rigidity of $R$,
$$\RHom_{A^\e}(A,R\otimes R^{\circ})\otimes_{A^\e} Q^\e
\cong R\otimes_{A^\e}Q^\e\cong X[n].$$
Combining these we have
$$\RHom_{Q^\e}(Q,X\otimes X^{\circ}[2n])\cong 
\RHom_{A^\e}(A,R\otimes R^{\circ})\otimes_{A^\e} Q^\e
\cong X[n].$$
After a complex shift we have 
$$\RHom_{Q^\e}(Q,X\otimes X^{\circ})\cong X[-n].$$
Since $X\cong Q^\tau$ and $X^{\circ}\cong (Q^\tau)^{\circ}$ as 
$Q$-bimodules, we have $\RHom_{Q^\e}(Q,Q^\e)\cong Q^\sigma[-n]$
where $\sigma=\tau^{-1}$. Hence $Q$ is rigid and $\fhtr Q=n$. 
The first assertion follows.

The proof of the graded case is similar.

In the graded case let $\tilde{R}=Q_{\rm gr}\otimes_A {R} \otimes_{A} 
Q_{\rm gr}$ where $Q_{\rm gr}=Q_{\rm gr}(A)$. By Lemma 4.1, $\tilde{R}$ is a 
graded dualizing complex over $Q_{\rm gr}$. Since $Q_{\rm gr}$ 
is noetherian and
has global dimension at most 1, $\tilde{R}$ is also an ungraded dualizing 
complex over $Q_{\rm gr}$. As said before, $\tilde{R}\cong 
Q_{\rm gr}^\sigma(l)[n]$; so $\tilde{R}$ is evenly localizable to $Q$. 
By hypothesis the $A^\e$-module $A$ is pseudo-coherent. Since $Q_{\rm gr}
\cong Q^\e_{\rm gr}\otimes_{A^\e}A$, $Q_{\rm gr}$ is pseudo-coherent over 
$Q^\e_{\rm gr}$. Hence we can apply the first assertion to $\tilde{R}$. 
The last assertion follows by the fact $\xi(\tilde{R})=
\xi_{\gr}(\tilde{R})=\xi_{\gr}(R)=\xi(R)$.
\end{proof}

To use Proposition \ref{xxprop4.2} we need to check the following:

\bigskip

(C1) the $A^\e$-module $A$ is pseudo-coherent;

(C2) there exists a rigid dualizing complex $R$ over $A$;

(C3) $R$ is evenly localizable to $Q$;

(C4) $\xi(R)$ is computable.

\bigskip

In the rest of this section we discuss (C1), (C2) and (C3).
First we consider condition (C1). 

If $A^\e$ is noetherian, then $A$ has a free resolution over 
$A^\e$ with each term being a finite free $A^\e$-module. So $A$ is 
pseudo-coherent over $A^\e$.

Let $A$ be a connected graded ring. Following \cite{VdB}, 
$A$ is called {\it $\Ext$-finite} if $\Ext^i_A(k,k)$ is 
finite dimensional over $k$ for all $i$. If $A$ is 
noetherian, then it is $\Ext$-finite. There are many
non-noetherian graded rings which are $\Ext$-finite. For 
example, if $A$ is AS regular (not necessarily noetherian), 
then $A$ is Ext-finite \cite[Proposition 3.1(3)]{SteZ}.

Let $F\to k$ be the minimal free resolution of $k$ as $A$-module. 
Then $F^{-i}\cong A\otimes  V_i$ where $V_i$ is the graded
vector space $\Tor^A_i(k,k)$. Hence $A$ is $\Ext$-finite
if and only if $\Tor^A_i(k,k)$ is finite dimensional for 
all $i$, if and only if $k$ is pseudo-coherent over $A$.

\begin{lemma}
\label{xxlem4.3} 
Let $A$ be a connected graded algebra.
Let $V_i=\Tor^A_i(k,k)$. 
\begin{enumerate}
\item 
The graded $A^{\e}$-module $A$ has a minimal graded free resolution 
$P$ such that $P^{-i}\cong A^{\e}\otimes  V_i$.
\item 
The projective dimension $\pdim_{A^\e} A$ is equal to the 
projective dimension $\pdim_A k(=\gldim A)$.
\item
If $A$ is $\Ext$-finite, then the $A^\e$-module $A$ is pseudo-coherent.
\item
Let $t$ be a homogeneous regular normal element of $A$. Let
$M$ be a graded pseudo-coherent $A/(t)$-module. Then $M$ is 
pseudo-coherent viewed as an $A$-module.
\item 
Let $t$ be a homogeneous regular normal element of $A$.
Then $\pdim_A k=\pdim_{A/(t)} k+1$ if $\pdim_{A/(t)} k$
is finite.
\end{enumerate}
\end{lemma}

We would like to remark that Lemma \ref{xxlem4.3}(b) is similar to
a result of Rouquier \cite[Lemma 7.2]{Rou} which says that 
$\pdim_{A^\e} A=\gldim A$ for finite dimensional algebras $A$ or 
commutative algebras $A$ essentially of finite type.

\begin{proof}[Proof of Lemma \ref{xxlem4.3}] 
(a) It suffices to show that $\Tor^A_i(k,k)\cong
\Tor^{A^\e}_i(k,A)$ as graded $k$-modules. Let $P$ be a minimal
graded free resolution of the graded $A^\e$-module $A$.
We think of $P$ as an $A$-bimodule free resolution of $A$. Restricted 
to the left (and to the right), the exact complex $P\to A\to 0$ is 
a split sequence since every term of it is a free $A$-module 
(respectively, free $A^\circ$-module). This implies that $k\otimes_{A^\circ}
(P\to A\to 0)$ is exact and hence $P':=k\otimes_{A^{\circ}} P$ 
is a free resolution of the $A$-module $k$. Hence
$$\Tor^A_i(k,k)=\HB^{-i}(k\otimes_A P')\cong 
\HB^{-i}(k\otimes_A (k\otimes_{A^{\circ}}P))=:(*).$$
Since $k\otimes_{A^{\circ}} P=(A\otimes  k)\otimes_{A^\e} P$
and $k\cong k\otimes_A(A\otimes k)$, we have
$$(*)\cong \HB^{-i}(k\otimes_A (A\otimes  k)\otimes_{A^\e} P)\cong
\HB^{-i}(k\otimes_{A^\e} P)=\Tor^{A^\e}(k,A).$$

(b,c) are immediate consequences of part (a).

(d,e) We use the double $\Tor$  spectral sequence:
\begin{equation}
\label{E4.3.1}
\Tor^B_p(\Tor^A_q(k,B), M)\Rightarrow_{p} \Tor^A_{n}(k,M)
\tag{E4.3.1}
\end{equation}
where $B=A/(t)$. Note that $\Tor^A_i(k,B)\cong k$ for $i=0,1$ and 
$\Tor^A_i(k,B)=0$ for all $i>1$. If $M$ is pseudo-coherent as $B$-module,
then $\Tor^B_i(k,M)$ is finite for all $i$. By the above spectral 
sequence $\Tor^A_i(k,M)$ is finite for all $i$. Thus we proved 
part (d).

For part (e) we assume that $\pdim_{B} k<\infty$ and let $M=k$ in the 
spectral sequence. Then we see that $\pdim_A k\leq 
\pdim_B k+1$. Also the spectral sequence \eqref{E4.3.1} implies that
$$\Tor^A_{n+1}(k,k)=\Tor^B_n(\Tor^A_1(k,B),k)\neq 0$$
for $n=\pdim_B k$. The assertion follows.
\end{proof}

Similarly, $k\otimes_A P$ is a minimal free resolution of $k$ 
as $A^{\circ}$-module. Combining Lemma \ref{xxlem4.3} with the comments 
before Lemma \ref{xxlem4.3} we have proved the following. 

\begin{proposition}
\label{xxprop4.4}
Condition (C1) holds if $A$ satisfies one of the following conditions:
\begin{enumerate}
\item 
$A^\e$ is noetherian.
\item
$A$ is connected graded and noetherian.
\item
$A$ is connected graded and AS regular.
\end{enumerate} 
\end{proposition}

Secondly we consider condition (C2). 

There are already many 
results about the existence of rigid dualizing complexes in
\cite{VdB,Ye1,YZ1}. For example, if $A$ has a filtration such 
that $\gr A$ is connected graded, noetherian, and AS 
Gorenstein, then $R={A^\sigma}[n]$ is a rigid dualizing complex
over $A$ where $n$ is the injective dimension of $\gr A$ 
\cite[Proposition 6.18(2)]{YZ1}. For non-noetherian AS regular algebras
rigid dualizing complexes also exist. 

\begin{proposition}
\label{xxprop4.5} Let $A$ be an AS regular 
algebra of (graded) global dimension $n$. 
\begin{enumerate}
\item
$A^\e$ is AS regular of global dimension $2n$.
\item
$\RHom_{A^\e}(A,A^\e)\cong {A^\tau}(-l)[-n]$ for some 
automorphism $\tau$, where $l$ is the number given in 
\textup{Definition \ref{xxdefn3.5}(c)}.
\item 
Let $R={A^{\sigma}}(l)[n]$ where $\tau=\sigma^{-1}$.
Then $R$ is a rigid dualizing complex over $A$.
\item
If $A$ is (graded) Goldie prime, then $Q(A)$ is rigid
and smooth and 
$$\fhtr A=\fhtr Q(A)=\fhtr_{\gr} Q_{\rm gr}(A)=n.$$
\end{enumerate}
\end{proposition}

\begin{proof} (a) For any  connected graded ring $B$,
$\gldim B=\pdim_B k$. It is clear that 
$$\pdim_{A^\e} k=\pdim_{A\otimes A^\circ} {_Ak}\otimes {_{A^\circ}k}
\leq \pdim_A k+\pdim_{A^\e} k=2n.$$
Hence $\gldim A^\e\leq 2n$. Note that $_Ak$ and $_{A^\circ}k$ 
have finite minimal free resolutions over $A$ and $A^\circ$
respectively \cite[Proposition 3.1(3)]{SteZ}.
The K{\"u}nneth formula and the AS Gorenstein property for 
$A$ and $A^\e$ imply that
$$\RHom_{A^\e}(k,A^\e)=\RHom_A(k,A)\otimes
\RHom_{A^\circ}(k,A^\circ)=k(-2l)[2n].$$
Therefore $A^\e$ is an AS regular of global dimension $2n$.

(b) Since $A$ is AS regular, the $A^\e$-module $A$ is 
pseudo-coherent [Proposition \ref{xxprop4.4}(c)]. Since $A$ has 
finite global dimension, $k$ has finite projective dimension. 
By Lemma \ref{xxlem2.2}(a),
$$\RHom_{A^\e}(A,A^\e)\otimes^{\rm L}_A k\cong 
\RHom_{A^\e}(A,A^\e\otimes^{\rm L}_A k)\cong
\RHom_{A^\e}(A,k\otimes  A^\circ)=:(*).$$
Since $A$ is AS regular, we have 
$$\RHom_A(k,A)=k(-l)[-n]\qquad {\text{or}}\qquad
\RHom_A(k,A(l)[n])=k,$$
where $(l)$ is the degree shift. Hence
$$k\otimes  A^\circ=\RHom_{A^\e}(k\otimes  A^\circ,
A\otimes  A^\circ(l)[n]).$$
The computation continues 
$$(*)\cong\RHom_{A^\e}(A,\RHom_{A^\e}(k\otimes  A^\circ,
A^\e(l)[n])\cong 
\RHom_{A^\e}((k\otimes  A^\circ)\otimes^{\rm L}_{A^\e}A, A^\e(l)[n])$$
The computation in the proof of Lemma \ref{xxlem4.3}(a) shows 
that
$$(k\otimes  A^\circ)\otimes^{\rm L}_{A^\e}A\cong
(k\otimes  A^\circ)\otimes_{A^e} P\cong 
k\otimes_A P \cong k.$$
Therefore
$$(*)\cong \RHom_{A^\e}((k\otimes  A^\circ)\otimes^{\rm L}_{A^\e}A, 
A^\e(l)[n])\cong \RHom_{A^\e}(k,A^\e(l)[n])=:(**).$$
By the proof of part (a), $A^\e$ is AS regular and 
$\RHom_{A^\e}(k,A^\e)=k(-2l)[-2n]$. Therefore 
$$(**)\cong k(-l)[-n].$$ Thus we have proved that
$$\RHom_{A^\e}(A,A^\e)\otimes^{\rm L}_A k\cong k(-l)[-n].$$

Since the free resolution $P$ of the $A^\e$-module $A$ is bounded with
each term being finite, $\RHom_{A^\e}(A,A^\e)\cong \Hom_{A^\e}(P,A^\e)
:=P^{\vee}$, which  is a bounded complex of finite free  
right $A^\e$-modules. Let $V$ be the minimal graded free
resolution of $P^{\vee}$ viewed as $A^{\circ}$-module complex. 
Note that the existence of $V$ follows from the facts that each term 
of $P^\vee$ is locally finite and that $(P^\vee)_{\ll 0}=0$. Then 
$$V\otimes^{\rm L}_{A}k\cong \RHom_{A^\e}(A,A^\e)
\otimes^{\rm L}_A k\cong k(-l)[-n]$$
which implies that $V\cong A^{\circ}(-l)[-n]$; or equivalently, 
$\RHom_{A^\e}(A,A^\e)\cong A^{\circ}(-l)[-n]$ as $A^{\circ}$-modules.
Similarly, $\RHom_{A^\e}(A,A^\e)\cong A(-l)[-n]$ as $A$-modules.
Thus 
$$\Ext^n_{A^\e}(A,A^\e)\cong {A^\sigma}(-l)$$
for some automorphism $\sigma$; and 
$$\Ext^i_{A^\e}(A,A^\e)=0$$ 
for all $i\neq n$. 

(c) The rigidifying isomorphism follows from part (b). It is easy to see 
that $R$ is a dualizing complex.

(d) Follows from part (c) and Proposition \ref{xxprop4.2}.
\end{proof}

Thirdly we look at condition (C3).

If $R={A^\sigma}[n]$, then $R$ is clearly evenly localizable 
to $Q(A)$ and $\xi(R)=n$. So we can now use Proposition 
\ref{xxprop4.2} to compute the $\fhtr$ (similar to Proposition 
\ref{xxprop4.5}(d)). 

Condition (C3) may hold even for $R\neq {A^\sigma}[n]$. 
In fact we don't have any examples of rigid dualizing complexes 
such that (C3) fails. Let's consider two cases: (i) $R$ is an 
Auslander dualizing complexes; and (ii) $A$ is a graded (or 
filtered) ring. We have already seen that Auslander dualizing 
complexes are evenly localizable to $Q(A)$ 
\cite[Proposition 3.3]{YZ4}. 

\begin{proposition}
\label{xxprop4.6} 
Let $A$ be a noetherian prime ring and let $Q=Q(A)$. Let $R$ be an
Auslander dualizing complex over $A$.
\begin{enumerate}
\item 
\cite[Proposition 3.3]{YZ4}
$R$ is evenly localizable to $Q$ and $\xi(R)=\Cdim A$.
\item
If $R$ is rigid and if $A^\e$-module $A$ is pseudo-coherent, then 
$Q$ is rigid and $\fhtr Q =\Cdim A$.
If $A$ is also graded, then $\fhtr Q_{\rm gr}(A)=\Cdim A$.
\end{enumerate}
\end{proposition} 

\begin{proof} (a) is \cite[Proposition 3.3]{YZ4}. To prove part 
(b) we use Proposition \ref{xxprop4.2}.
\end{proof}

A useful consequence is the following.  Recall from \cite{YZ1} 
that a noetherian connected graded ring $A$ is said to {\it 
have enough normal elements} if every non-simple prime graded 
factor ring $A/I$ contains a nonzero normal element of positive 
degree. 

\begin{corollary} 
\label{xxcor4.7}
Let $A$ be a prime ring with a noetherian connected filtration 
such that $\gr A$ has enough normal elements. Then $Q$ is rigid and 
$\fhtr Q=\GKdim A$.
\end{corollary}

\begin{proof} By \cite[Section 4]{ASZ}, $A\otimes B$ is noetherian
for any noetherian algebra $B$. So $A^\e$ is noetherian and hence
$A$ is pseudo-coherent. By \cite[Corollary 6.9]{YZ1}, $A$ has an Auslander
rigid dualizing complex $R$ such that $\Cdim=\GKdim$. The assertion
follows from Proposition 4.6(b).
\end{proof} 

Finally we mention that (C3) holds for graded noetherian rings.

\begin{lemma}
\label{xxlem4.8} Let $A$ be a graded noetherian prime ring with a 
balanced dualizing complex $R$. Then $R$ is evenly localizable
to $Q(A)$ and to $Q_{\rm gr}(A)$.
\end{lemma}

The proof follows from Lemma \ref{xxlem6.3}(b,c) below where the 
filtered case is considered.

\section{Proof of Theorem \ref{xxthm0.2}}
\label{xxsec5}

In this section we prove Theorem \ref{xxthm0.2}. Let $U$ denote 
some division $k$-algebra. We sometimes use $A_U$ to denote the ring 
$A\otimes  U$. 

\begin{lemma}
\label{xxlem5.1} Let $A$ be a connected graded ring and let $U$ 
be a division ring. Then the ungraded global dimension $\gldim 
A_U$ is equal to the graded projective dimension 
$\pdim_{A} k$.
\end{lemma}

\begin{proof} It is clear that 
$$\pdim_A k\leq \gldim A\leq \gldim A_U.$$
So it suffices to show that $\pdim_A k\geq \gldim A_U$.
Hence we may assume that $\pdim_A k<\infty$. Let $P$ be a
minimal free $A^\e$-resolution of $A$. By Lemma \ref{xxlem4.3}(b),
the length of $P$ is equal to $\pdim_{A} k$. Let $P_U:=P\otimes  U$
and consider $U$ as a $U^\e$-module. Then $P_U$ is a complex of 
$A_U^\e$-module. Let $M$ be any ungraded $A_U$-module, we claim that
$M\otimes_{A^{\circ}_U} P_U$ is a free $A_U$-module resolution of 
$M$. By the reasoning given in the proof of Lemma \ref{xxlem4.3}(a),
$P\to A\to 0$ is a split exact sequence of $A^\circ$-modules. 
Hence $P_U\to A_U\to 0$ is a split exact sequence of 
$A^{\circ}_U$-modules. Thus $M\otimes_{A^{\circ}_U} P_U\to M\to 0$ 
is exact, or $M\otimes_{A^{\circ}_U} P_U$ is a resolution of $M$. 
Now every term in $M\otimes_{A^{\circ}_U} P_U$ is a direct
sum of copies of
$$M\otimes_{A^{\circ}_U}(A^\e\otimes  U)\cong M\otimes  A$$
where $M$ is viewed as a $U$-module. Since $U$ is a division ring, 
$M$ is free over $U$, and whence $M\otimes  A$ is free over $A_U$. 
Therefore the $M\otimes_{A^{\circ}_U} P_U$ is a projective 
resolution of $M$. Consequently, the projective dimension of $M$ 
is bounded by the length of $P$. Thus $\gldim A_U\leq \pdim_{A^\e}A=
\pdim_A k$. 
\end{proof}

\begin{proposition}
\label{xxprop5.2}
If $A$ is a connected graded Goldie prime ring of finite global
dimension, then $\shtr Q(A)\leq \gldim A$.
\end{proposition}

\begin{proof} The assertion follows from the (in)equalities
$$\injdim Q(A)\otimes U\leq \gldim Q(A)\otimes U \leq \gldim 
A\otimes U=\gldim A$$
where the last equality is Lemma \ref{xxlem5.1}.
\end{proof}

We now prove a restatement of Theorem \ref{xxthm0.2}.

\begin{theorem}
\label{xxthm5.3} Let $A$ be an AS regular graded ring that is 
Goldie prime. Then $Q(A)$ is rigid, smooth and homologically 
uniform and $\htr Q(A)=\gldim A$.
\end{theorem}

\begin{proof} Since $A$ is AS regular, so is $A\otimes A^\circ$
[Proposition \ref{xxprop4.5}(a)].
By Lemma \ref{xxlem1.4}(d), $Q\otimes Q^{\circ}$ has
finite global dimension where $Q:=Q(A)$. Hence $Q$ is 
is smooth [Lemma \ref{xxlem1.3}(b)]. By 
Proposition \ref{xxprop4.5}(d), $Q$ is rigid and $\fhtr Q=\gldim A$. 
By Proposition \ref{xxprop5.2}, $\shtr Q\leq \gldim A$. 
The assertion follows from Lemma \ref{xxlem1.2}(a).
\end{proof}

We do not assume that $A$ is noetherian in the above statement. 
There is no reason to expect that an AS regular Goldie prime ring
(or even Ore domain) should be noetherian or have finite GK-dimension. 
There are many connected graded Ore domains have exponential 
growth (hence not noetherian). For example any Rees ring of an 
affine Ore domain of exponential growth is connected graded with 
exponential growth. 

Next we consider the degree zero part of the graded quotient ring.
Let $A$ be a graded prime ring and let $Q_{\rm gr}(A)$ be the graded
Goldie quotient ring. Let $Q_0(A)$ (or $Q_0$ for short) be the degree
zero part of $Q_{\rm gr}(A)$. In general,
$Q_0$ is semisimple artinian and it is a finite direct sum of
simple artinian rings which are isomorphic to each other. We define
$\htr Q_0$ to be the $\htr$ of one copy of its simple artinian 
summands. If $A$ is an Ore domain, then $Q_0$ is a division algebra. 
We now prove a version of Theorem \ref{xxthm5.3} for $Q_0$.

\begin{proposition}
\label{xxprop5.4} 
Let $A$ be an AS regular graded ring that is Goldie prime. 
Suppose that $Q_{\rm gr}(A)=Q_0[t^{\pm 1};\alpha]$ for some 
$t\in A_{\geq 1}$ and some automorphism $\alpha$. If $Q_0$ is 
doubly noetherian, then $Q_0$ is rigid, smooth and homologically 
uniform and $\htr Q_0=\gldim A-1$.
\end{proposition}

\begin{proof} Let $Q$ be the Goldie quotient ring of $A$. Then
$Q=Q_0(t;\alpha)$. By Theorem \ref{xxthm5.3}, $Q$ is rigid,
smooth and homologically uniform. By Lemma \ref{xxlem2.1}(c), 
$Q_0$ is smooth; and by hypothesis, $Q_0$ is doubly noetherian. 
Lemma \ref{xxlem1.3} implies that $Q_0$ is homologically uniform. 
To show $\htr Q_0=\gldim A-1$, it suffices to show that 
$\fhtr Q_0=\fhtr Q-1$, which follows from Proposition 
\ref{xxprop2.6}. Also by Proposition \ref{xxprop2.6}, $Q_0$ is 
rigid. 
\end{proof}

\section{Proof of Theorem \ref{xxthm0.3}}
\label{xxsec6}

At the end of this section we prove Theorem \ref{xxthm0.3}. In 
order to prove Theorem \ref{xxthm0.3} we need to review (and extend) 
some work of Van den Bergh \cite{VdB} and of the first author \cite{Ye1} 
on local duality for graded modules, and to study filtered rings. 

When applied to graded modules $\Hom$ and $\otimes$ and their 
derived functors will be in the graded sense. If $M$ is a graded 
$k$-vector space, then $M'$ is the graded $k$-linear dual of $M$. 

Recall that a connected graded ring $A$ is $\Ext$-finite if 
$\Ext^i_A(k,k)$ is finite for all $i$. A consequence of this condition 
is that $A/A_{\geq n}$ is pseudo-coherent over $A$ for every $n\geq 1$. 
Let $U$ be any division algebra over $k$. 
By tensoring with $U$  we see that $A_U/(A_U)_{\geq n}$ 
pseudo-coherent over $A_U$ for all $n$. The trivial graded 
$A_U$-module $A_U/(A_U)_{\geq 1}$ is also denoted by $U$.

\begin{definition}
\label{xxdefn6.1}
Let $B$ be any $\mathbb N$-graded ring (not necessarily connected 
graded) and let $\fm=B_{\geq 1}$. 
\begin{enumerate}
\item
For any graded $B$-module $M$, the {\it ${\fm}$-torsion functor} 
$\Gamma_{\fm}$ is defined to be
$$\Gamma_{\fm} (M)=\{x\in M\; |\; {\fm}^n x=0, \;{\text{for}}\; n\gg 0\}.$$
The right derived functor of $\Gamma_{\fm}$, denoted by $\R\Gamma_{\fm}$,
is defined on the derived category $\D^{+}(\GrMod B)$. 
\item
The {\it $i$th local cohomology} of $X\in \D^+(\GrMod B)$ 
is defined to be
$$\HB_{\fm}^i(X)=\R^i \Gamma_{\fm}(X).$$
\item
The {\it local cohomological dimension} of a graded $B$-module $M$ is
defined to be
$$\lcd M=\sup\{i\; |\; \HB_{\fm}^i(M)\neq 0\}.$$
\item
The {\it cohomological dimension} of $\Gamma_{\fm}$, also called 
the {\it cohomological dimension} of $B$, is defined to be
$$\cd \Gamma_{\fm}=\cd B=
\sup\{\lcd M \; |\; {\text{for all graded $B$-modules $M$}}\}.$$
\end{enumerate}
\end{definition}

Obviously, $\Gamma_{\fm}(M)=\varinjlim \Hom_B(B/{\fm}^n, M)$, which 
implies that 
$$\HB_{\fm}^i(X)=\varinjlim \Ext^i_B(B/{\fm}^n, X)$$ 
for all 
$X\in \D^{+}(\GrMod B)$. If $B$ is left noetherian and $\cd \Gamma_{\fm}
<\infty$, then 
$$\cd \Gamma_{\fm} =\lcd B=\sup\{i\;|\; \HB_{\fm}^i(B)\neq 0\}.$$

\begin{lemma}
\label{xxlem6.2} Assume $A$ is $\Ext$-finite. 
Let $R=\operatorname{R\Gamma_{\fm}}(A)'$ and suppose it is locally 
finite. Let $R_U=R\otimes  U$.
\begin{enumerate}
\item
Let $E=A'\otimes  U$. Then $E$ is the graded 
injective hull of the graded trivial 
$A_U$-module $U$ and $\Hom_U(M,U)\cong \Hom_{A_U}(M,E)$ for all graded
$A_U$-modules $M$.
\item
Let $M$ be a graded $A_U$-module. Then 
$\operatorname{R\Gamma_{\fm_U}}(M)\cong \operatorname{R\Gamma_{\fm}}(M)$.
As a consequence, $\cd A=\cd A_U$.
\item 
Suppose $\cd A<\infty$. Then 
$$\Hom_U(\operatorname{R\Gamma_{\fm_U}}(M),U)\cong \RHom_{A_U}(M,R_U)$$
for $M\in \D^{b}(\GrMod A)$. 
\item 
Suppose $\cd A<\infty$. Then $R_U$ has 
graded injective dimension 0.
\item
If $A$ is noetherian with balanced dualizing complex $R$, then $R_U$ is 
a graded dualizing complex over $A_U$.
\item
If $A$ is AS Gorenstein, then the graded injective dimension
of $A_U$ is equal to the graded injective dimension of $A$.
\end{enumerate}
\end{lemma}

\begin{proof} (a) Consider the exact functor $F: M\mapsto \Hom_U(M,U)$
from $\GrMod A_U$ to $\GrMod A_U^\circ$. This functor sends coproducts to 
products. By Watts' theorem \cite[Theorem 3.36]{Rot}, $F$ is 
equivalent to $M\mapsto \Hom_{A_U}(M,E_U)$ where 
$$E_U=\Hom_U(A_U,U)\cong A'\otimes  U=E.$$ 
Since the functor is exact, $E_U$ is injective. The socle of $E_U$ is 
the trivial module $U$, so $E_U$ is the injective hull of $U$.

(b) Let $M$ be a graded $A_U$-module. Then the $\Hom$-$\otimes$
adjunction implies that 
$$\Ext^i_{A_U}(A_U/(A_U)_{\geq n},M)=\Ext^i_{A}(A/A_{\geq n},M).$$
Therefore 
$$\operatorname{R\Gamma_{\fm_U}}(M)=\operatorname{R\Gamma_{\fm}}(M).$$
The assertion follows.

(c) Since $A$ is Ext-finite, $\operatorname{R\Gamma_{\fm}}(-)$ 
commutes with coproducts \cite[Lemma 4.3]{VdB}. By part (b) 
$\operatorname{R\Gamma_{\fm_U}}(-)$ commutes with coproducts. 

Using part (b) and the fact that $R$ is locally finite, we have
$$\Hom_{U}(\operatorname{R\Gamma_{\fm_U}}(A_U),U)\cong
\Hom_{U}(\operatorname{R\Gamma_{\fm}}(A_U),U)\cong
\Hom_{U}(R'\otimes U,U)\cong R\otimes  U=R_U.$$
Let $F$ be a bounded resolution of $A_U$ as graded 
$A^\e_U$-module whose restriction consists of $\Gamma_{\fm_U}$-acyclic
$A$-modules. Then we have $\operatorname{R\Gamma_{\fm_U}}(A)=
\Gamma_{\fm_U}(F)$. Let $K$ be a projective resolution of 
$M$. Since $\operatorname{R\Gamma_{\fm_U}}(-)$ commutes with 
coproducts, $F\otimes_{A_U} K$ is $\Gamma_{\fm_U}$-acyclic. This 
implies that 
$$\operatorname{R\Gamma_{\fm_U}}(M)=
\Gamma_{\fm_U}(F\otimes_{A_U} K)=\Gamma_{\fm_U}(F)\otimes_{A_U} K.$$
Let $E$ be the injective resolution of $U$ as a graded $A_U$-module.
By part (a) we have 
$$\RHom_{U}(\operatorname{R\Gamma_{\fm_U}}(M),U)
\cong \Hom_{A_U}(\Gamma_{\fm_U}(F)\otimes_{A_U} K,E)=:(*).$$
By the $\Hom$-$\otimes$ adjunction we have 
$$(*)\cong \Hom_{A_U}(K,\Hom_{A_U}(\Gamma_{\fm_U}(F), E))
\cong \RHom_{A_U}(K, R_U)\cong\RHom_{A_U}(M, R_U)$$
where the middle isomorphism follows from 
$$\Hom_{A_U}(\Gamma_{\fm_U}(F), E)\cong \Hom_{A_U}(
\operatorname{R\Gamma}_{\fm_U}(A_U), E)
\cong R_U.$$

(d) Since the complex $\operatorname{R\Gamma_{\fm_U}}(M)$ lives in 
non-negative positions, the ``dual'' complex 
$\Hom_{U}(\operatorname{R\Gamma_{\fm_U}}(M),U)$ 
lives in the non-positive positions. Hence by part (c), $R_U$ has 
injective dimension less than or equal to 0. If $M$ is $\fm$-torsion, 
then, by part (c), $\Ext^0_{A_U}(M,R_U)\neq 0$. The assertion follows.

Note that if $M$ is $\fm$-torsion-free then $\Ext^0_{A_U}(M,R_U)=0$. 

(e) By definition $R$ is pseudo-coherent over $A$ on both sides. 
So $R_U$ is pseudo-coherent over $A_U$ on both sides. By part (d) 
$R_U$ has finite injective dimension on the left; by symmetry also
on the right. Finally $\RHom_{A_U}(R_U,R_U)\cong A_U$ follows from 
the fact $R$ is pseudo-coherent and $R$ is a dualizing complex over $A$.
By symmetry, $\RHom_{A^\circ_U}(R_U,R_U)\cong A_U$. Therefore
$R_U$ is a dualizing complex over $A_U$.

(f) If $A$ is AS Gorenstein (and $\Ext$-finite),
then the balanced dualizing complex over $A$ is $R=A^\sigma(l)[n]$ 
where $n$ is the graded injective dimension of $A$. The assertion 
follows from part (d).
\end{proof}

An ascending $\mathbb N$-filtration 
$F=\{F_i A\}_{i\geq 0}$ on a ring $A$ is called {\it connected filtration}
(respectively, {\it noetherian connected filtration})
if 
\begin{enumerate}
\item[(i)] $1\in F_0 A$, 
\item[(ii)] $F_i A\; F_j A \subset F_{i+j} A$, 
\item[(iii)] 
$A=\bcup_{i\geq 0} F_i A$, and 
\item[(iv)] the associated graded ring 
$$\gr A:=\bigoplus_{i=0}^{\infty} F_i A/F_{i-1} A$$ 
is connected graded (respectively, connected graded and noetherian).
\end{enumerate}
The Rees ring of $A$ with a given filtration $F$ defined to be
$$L:=\bigoplus_{i=0}^{\infty} (F_i A)t^i.$$ 
So $L$ is a subring of $A[t]$ such that $L/(t)=\gr A$ and $L/(t-1)=A$. 
By \cite[Theorem 8.2]{ATV}, $\gr A$ is noetherian if and only 
if $L$ is.

An $A$-bimodule is called {\it filtered finite} if there 
is a filtration on $M$ compatible with the filtration on $A$ such 
that $\gr M$ is a finite left and a finite right graded $\gr A$-module.
If $A$ is connected graded, then $A$ has an obvious filtration such 
that $\gr A=A$. Sometimes we view a graded ring $A$ as a filtered 
ring so that we can pass some graded properties to the ungraded setting. 

\begin{lemma}
\label{xxlem6.3} Let $A$ be a ring with a connected filtration and 
let $L$ be the Rees ring. In parts (c-f) assume that $A$ has a 
noetherian connected filtration and that $A$ is prime.
Let $R$ be a rigid dualizing complex over $A$ (if it exists). 
\begin{enumerate}
\item
If $\gr A$ is $\Ext$-finite, then so is $L$.
\item
If $A$ is Goldie prime, then $L$ is both graded Goldie prime
and ungraded Goldie prime. Furthermore $Q(L)=Q(A)(t)$. 
\item
Every filtered finite bimodule is
evenly localizable.
\item
If $\gr A$ has a balanced dualizing complex, 
then any rigid dualizing complex $R$ over $A$ is 
evenly localizable.
\item
If $\gr A$ has a balanced dualizing complex, then
$$0\leq \xi(R)\leq \cd \gr A=-\inf\{i \;|\; \HB^i(R)\neq 0\}.$$
\item
If $\gr A$ is noetherian and AS Gorenstein, then $\xi(R)=\injdim
\gr A$. 
\end{enumerate}
\end{lemma}

\begin{proof} (a) We know that $t$ is a central regular element in $L$
such that $L/(t)\cong \gr A$. Applying  Lemma \ref{xxlem4.3}(d)
to the trivial module $k$ the assertion follows.

(b) For every regular element $x\in A$, $x t^i$ 
for some $i\geq 0$ is a homogeneous regular element in $L$. The 
set of homogeneous regular elements in $L$ form an Ore set. By 
inverting all homogeneous elements we obtain that
$$Q_{\rm gr}(L)=Q(A)[t^{\pm 1}],$$ 
which is prime and graded simple artinian. So $L$ is both 
graded Goldie prime and ungraded Goldie prime, and $Q(L)=Q(A)(t)$. 

(c) Let $M$ be a filtered finite $A$-bimodule. 
By \cite[Lemma 3.1]{SZ}, $M$ is left Goldie torsion if and only if
$M$ is a right Goldie torsion. By passing to the factor module 
of $M$ modulo the largest Goldie torsion $A$-submodule we may assume 
that $M$ is two-sided Goldie torsionfree. Since $Q\otimes_A M$ 
is an artinian left $Q$-module 
and Goldie torsionfree right $A$-module, any regular element 
in $A^\circ$ acts on $Q\otimes_A M$ bijectively. Hence 
$Q\otimes_A M\cong Q\otimes_A M\otimes_A Q$.
Similarly, $M\otimes_A Q\cong Q\otimes_A M\otimes_A Q$. 
Thus $M$ is evenly localizable.

(d) It follows from the construction of the rigid dualizing 
complex $R$ in \cite[Theorem 6.2]{YZ1} that $\HB^i(R)$ is filtered 
finite. By part (c) each $\HB^i(R)$ is evenly localizable. So 
$R$ is evenly localizable.

(e) Let $R_{\gr A}$ be a balanced dualizing complex over $\gr A$.
By \cite[Theorem 6.3]{VdB} $R_{\gr A}=
\operatorname{R\Gamma_{\fm}}(\gr A)'$. Since 
$\Gamma_{\fm}$ has finite cohomological dimension,
$\cd \gr A= \cd \Gamma_{\fm}=
-\inf\{i \;|\; \HB^i(R_{\gr A})\neq 0\}$.
By the construction of $R$ in \cite[Theorem 6.2]{YZ1} one 
sees that 
$$\inf\{i \;|\; \HB^i(R)\neq 0\}=
\inf\{i \;|\; \HB^i(R_{\gr A})\neq 0\}=-\cd \gr A.$$ 
Since $R_{\gr A}$ has injective dimension $0$, $R$ has injective 
dimension at most 0 by the construction. If $\HB^i(R)\neq 0$ then 
$i$ lies in between $-(\cd \gr A)$ and $0$. The assertion follows.

(f) When $A$ is filtered AS Gorenstein, then $R={A^\sigma}[n]$
where $n=\injdim \gr A$ \cite[Proposition 6.18]{YZ1}. Hence 
$\xi(R)=\injdim \gr A$. 
\end{proof}

\begin{proposition}
\label{xxprop6.4} Let $A$ be a Goldie prime ring. Suppose $A$ has a 
connected filtration such that $\gr A$ is AS regular. If $Q(A)$ is 
doubly noetherian, then $Q(A)$ is smooth, rigid and homologically 
uniform, and $\htr Q(A)=\gldim \gr A$. 
\end{proposition}

\begin{proof} Let $L$ be the Rees ring of $A$. By Lemma 
\ref{xxlem4.3}(e) and the Rees lemma, $L$ is AS regular
and $\gldim L=\gldim \gr A+1$. By Lemma 
\ref{xxlem6.3}(b), $L$ is Goldie prime. Since $Q_{\rm gr}(L)=
Q(A)[t^{\pm 1}]$ and $Q(A)$ is doubly noetherian, the assertion follows
from Proposition \ref{xxprop5.4}.
\end{proof}

\begin{proposition}
\label{xxprop6.5} Let $A$ be a prime ring with a noetherian 
connected filtration such that $\gr A$ has a balanced dualizing 
complex. Let $R$ be a rigid dualizing complex over $A$. 
If $Q$ is doubly noetherian, then $\fhtr Q=\xi(R)$.
\end{proposition}

\begin{proof} Let $L$ be the Rees ring and let $R_L$ be a  
balanced (and rigid) dualizing complex over $L$. Since
$L$ is noetherian, it is pseudo-coherent as $L^\e$-module
(see Proposition \ref{xxprop4.4}(b)). All conditions 
in Proposition \ref{xxprop4.2} holds for $L$ in both the graded 
and the ungraded settings. Hence $\fhtr_{\gr} Q_{\rm gr}(L)=\xi_{\gr}(R_L)
=\fhtr Q(L)=\xi(R_L)$. 

By \cite[Theorem 6.2]{YZ1}, a rigid dualizing complex $R_A$ over $A$
is given by the degree zero part of $R_L[t^{-1}][-1]$. 
Thus $\HB^i(R_A)$ is not Goldie $A$-torsion if and only if 
$\HB^{i-1}(R_L)$ is not Goldie $L$-torsion. This implies that
$\xi(R_A)=\xi_{\gr}(R_L)-1$. It remains to show that
$\fhtr Q(A)=\fhtr Q_{\rm gr}(L)-1$. But this follows from 
Proposition \ref{xxprop2.6}.
\end{proof}

Next we study the injective dimension of $Q(A)\otimes U$.
The following lemma is \cite[Theorem 1.3]{GJ}.

\begin{lemma} 
\label{xxlem6.6}
\cite[Theorem 1.3]{GJ}.
Let $C$ be a noetherian ring and $t$ a central 
element of $C$. If $I$ is an injective $C$-module, then $I[t^{-1}]$
is an injective $C[t^{-1}]$-module. A graded version of the assertion
also holds.
\end{lemma}

\begin{proposition}
\label{xxprop6.7} Let $A$ be a filtered ring and let $L$ be 
the Rees ring. Assume $L$ is noetherian and has a balanced 
dualizing complex $R_L$. Let $U$ be a division ring such that 
$L_U$ is noetherian.
Let $R_A$ be a rigid dualizing complex over $A$. 
Then $R_A\otimes U$ is a dualizing complex over $A$ of injective
dimension at most 0.
\end{proposition}

\begin{proof} By \cite[Theorem 6.2]{YZ1}, 
$R_A=(R_{L}[t^{-1}])_0[-1]$. Hence, for all $U$,
$$R_A\otimes U=((R_L\otimes U)[t^{-1}])_0[-1].$$
Since $L_U$ is noetherian, $A_U$ is noetherian. This implies
that $R_A\otimes U$ is pseudo-coherent on both sides. 
Since $\RHom_A(R_A,R_A)\cong
A$, Lemma \ref{xxlem2.2} implies that 
$$\RHom_{A_U}(R_A\otimes U,R_A\otimes U)\cong \RHom_{A}(R_A,R_A)\otimes U
\cong A_U.$$ It remains to show that the injective
dimension of $R_A\otimes U$ is at most 0 on both sides. By Lemma 
\ref{xxlem6.2}(d), $R_L\otimes U$ has graded injective dimension
0. Let $I$ be the minimal graded injective dimension of $R_L\otimes U$.

We claim that $I^0\cong L'\otimes U$. Since $L\otimes U$ is 
noetherian, $I^0$ is a direct sum of indecomposable
injectives, say $\boplus J_i$. By Lemma \ref{xxlem6.2}(c),
$$\Ext^0_{L_U}(U,R_L\otimes U)=\Hom_U(U,U)=U.$$
This shows that $I^0$ contains only one copy of $L'\otimes U$.
If $I^0\neq L'\otimes U$, then there is an $\fm$-torsionfree graded
$A_U$-module $M$ such that $\Ext^0_{L_U}(M,R_L\otimes U)\neq 0$. 
But this contradicts Lemma \ref{xxlem6.2}(c) since 
$\Gamma_{\fm}(M)=0$. So we proved our claim.

Since $I^0$ is $t$-torsion, $I^0[t^{-1}]=0$; and by Lemma \ref{xxlem6.6} 
the complex $(R_L\otimes U)[t^{-1}]$ of $\mathbb Z$-graded 
$L[t^{-1}]$-modules has injective dimension at most $-1$. Since 
$L[t^{-1}]$ is strongly graded and $(L[t^{-1}])_0=A$, the complex 
$((R_L\otimes U)[t^{-1}])_0[-1]$ of $A$-modules has injective dimension 
at most 0.
\end{proof}

\begin{proposition}
\label{xxprop6.8} Let $A$ be a prime ring with a noetherian 
connected filtration such that $\gr A$ has a balanced dualizing 
complex. Let $U$ be a division ring such that $\gr A\otimes U$
is noetherian. Let $R$ be a rigid dualizing complex over $A$. 
Then $\injdim Q(A)\otimes U\leq \xi(R)$.
\end{proposition}

\begin{proof} It is clear that $\gr A\otimes U$ is noetherian
if and only if $L\otimes U$ is noetherian where $L$ is the Rees 
ring of $A$. In this case $A\otimes U$ is also noetherian.
By Proposition \ref{xxprop6.7}, $R_A\otimes U$ is a dualizing 
complex over $A_U$ with injective dimension $\leq 0$. 
By Lemma \ref{xxlem4.1}, $(Q\otimes_A R_A\otimes_A Q)\otimes U$ 
is a dualizing complex over $Q\otimes U$ where $Q=Q(A)$. 
Since $Q\otimes_A R_A\otimes_A Q\cong {Q^\sigma}[-d]$ where 
$d=\xi(R_A)$, $(Q\otimes_A R_A\otimes_A Q )\otimes U\cong 
{Q^\sigma}[-d] \otimes U$. The injective dimension will not 
increase under localization by Lemma \ref{xxlem2.3}. So 
$\injdim_{Q\otimes U} (Q\otimes U[-d])\leq 0$. The assertion 
follows by a complex shift. 
\end{proof}

We are now ready to prove a generalization of Theorem \ref{xxthm0.3}.

\begin{theorem}
\label{xxthm6.9}
Let $A$ be a filtered Goldie prime ring with noetherian connected 
filtration. Let $Q=Q(A)$.
\begin{enumerate}
\item 
If $\gr A$ has a enough normal elements, then $Q$ is rationally 
noetherian, rigid and 
homologically uniform, and $\htr Q =\GKdim A$.
\item
If $\gr A$ has an Auslander balanced dualizing complex
and $\gr A\otimes Q^\circ$ is noetherian, then $Q$ is rigid and 
$\fhtr Q=\htr Q =\cd \gr A$.
\item
If $\gr A$ is AS Gorenstein and $\gr A\otimes Q^\circ$ is noetherian, 
then $Q$ is rigid and $\fhtr Q=\htr Q=\injdim \gr A$.
\end{enumerate}
\end{theorem}

\begin{proof}
(a) If $\gr A$ has enough normal elements, then $\gr A\otimes U$
is noetherian for all $k$-algebras $U$ \cite[Proposition 4.9]{ASZ}. 
Consequently, $A$ and $Q$ are rationally noetherian. Also 
by \cite[Corollary 6.9]{YZ1}, $A$ has an Auslander, $\GKdim$-Macaulay, 
rigid dualizing complex $R_A$. Hence 
$$\xi(R_A)=\cd \gr A=\GKdim \gr A=\GKdim A.$$
By Proposition \ref{xxprop6.8}, $\injdim Q\otimes U\leq \xi(R_A)=
\GKdim A$ for all division rings $U$. This says that $\shtr Q\leq
\GKdim A$. The assertion follows from Lemma \ref{xxlem1.2} and 
Corollary \ref{xxcor4.7}.

(b) By \cite[Corollary 6.8]{YZ1} and the proof of \cite[Theorem 6.2]{YZ1}, 
$A$ has an Auslander dualizing complex $R_A$, and  $\cd \gr A=\Cdim \gr A
=\Cdim A$. By Proposition \ref{xxprop4.6}(a), $\xi(R_A)=\Cdim A$.  
By Proposition \ref{xxprop6.8}, $\injdim Q\otimes Q^\circ\leq \xi(R_A)$
since $\gr A \otimes Q^\circ$ is noetherian. By Proposition \ref{xxprop6.5},
$\fhtr Q=\xi(R_A)$. Hence $\fhtr Q=\htr Q=\xi(R)=\cd \gr A$. By Proposition
\ref{xxprop4.2} $Q$ is rigid.

(c) Similar to the proof of part (b). Using Proposition \ref{xxprop4.2},
Lemmas \ref{xxlem6.3}(f) and \ref{xxlem1.2}(a) and Proposition 
\ref{xxprop6.8}, we have 
$$\xi(R)=\fhtr Q\leq \htr Q=\injdim Q\otimes Q^\circ\leq
\xi(R)=\injdim \gr A.$$
Hence $\fhtr Q=\htr Q=\injdim \gr A$. By Proposition
\ref{xxprop4.2} $Q$ is rigid.
\end{proof}

Theorem \ref{xxthm0.3} is an immediate consequence of Theorem \ref{xxthm6.9}.
Theorem \ref{xxthm6.9} also applies to affine prime PI rings and various
quantum algebras which have noetherian filtrations
such that $\gr A$ has Auslander dualizing complexes. Here we
give an example of this kind.

\begin{example}
\label{xxex6.10} Let $Q$ be a division algebra finite over its center
$C$, and assume $C$ is finitely generated as a field. It is easy to pick  
an affine prime PI subalgebra $A$ of $Q$ such that
$Q$ is the quotient ring of $A$, and $A$ has a filtration such that
$\gr A$ is connected graded noetherian and affine PI. By Theorem 
\ref{xxthm6.9}(a), $Q$ is a rationally noetherian, rigid and 
homologically uniform; and $\htr Q=\GKdim Q=\tr Q$.
\end{example}

The following corollary follows from Theorem \ref{xxthm6.9}(a);
the proof is omitted.

\begin{corollary}
\label{xxcor6.11}
The simple artinian rings given in \textup{Example 
\ref{xxex1.9}(a,c,d,e,f,g)}
are rigid.
\end{corollary}

\section{Examples}
\label{xxsec7}

In this section we will give some examples to show that $\htr$
can be different from other versions of transcendence degrees
for certain division algebras.

\begin{proposition}
\label{xxprop7.1}
Let $F$ be a countably infinite dimensional separable algebraic 
field extension of $k$. 
\begin{enumerate}
\item
$F$ is smooth, homologically uniform and $\htr F=\gldim F^\e=1$.
\item
$F$ is not doubly noetherian and
$$\htr F =1>0=\sup\{\htr G_i \;|\; {\textrm{for all}}\; G_i\subset F
\;{\text{with}}\; \dim_k G_i<\infty\}.$$
\item
$\thtr F=\tr F=0$.
\item
$F$ is not rigid. 
\end{enumerate}
\end{proposition}

The proof of the above proposition follows from several lemmas below. 

\begin{lemma}
\label{xxlem7.2} Let $F$ be as in Proposition \ref{xxprop7.1}.
Then $F$ is smooth of $\shtr F\leq 1$ and $\gldim F^\e=\pdim_{F^\e} F=1$.
\end{lemma}

\begin{proof} By \cite[Theorem 10]{Os}, $\pdim_{F^\e} F=1$. 
The assertion follows from Lemma \ref{xxlem1.3}(a).
\end{proof}

Consider the short exact sequence
\begin{equation}
\label{E7.2.1}
0\to J\to F^\e\to F\to 0
\tag{E7.2.1}
\end{equation}
where the map $F^\e\to F$ is the multiplication and $J$ is 
the kernel of this map. 

The following Lemma \ref{xxlem7.3}(a) is due to Ken Goodearl.
The authors thank him for providing the result. 

\begin{lemma}
\label{xxlem7.3} Let $F$ be as in Proposition \ref{xxprop7.1}
and let $J$ be as in \eqref{E7.2.1}. Then 
\begin{enumerate}
\item 
There is an infinite sequence of nonzero orthogonal 
idempotents $\{e_i\}_{i=1}^{\infty}\subset F^\e$ such that 
$J=\boplus_{i=1}^{\infty} e_i F^\e$.
\item 
$\Hom_{F^\e}(F,J)=\Hom_{F^\e}(J,F)=0$.
\item
$\Ext^1_{F^\e}(F,F)=\Hom_{F^\e}(F,F^\e)=0$.
\item
$\Ext^1_{F^\e}(F,F^\e)$ is infinite dimensional over $F$.  
\end{enumerate}
\end{lemma}

\begin{proof} 
(a) Write $F$ as the union of countable sequence of finite 
dimensional subfields
$$k\subset F_1\subset F_2\subset \cdots \subset F.$$
Let $J_i$ be the kernel of the map $F^\e_i\to F_i$.
Since $F$ is separable, $F_i^\e$ is a finite direct sum of
field extensions of $k$. This implies that there are idempotents
$u_i$ and $v_i=1-u_i$ in $F^\e_i$ such that $F^\e_i=
u_i F^\e_i \boplus v_i F^\e_i$ where $J_i=u_i F^\e_i$ 
and $F_i \cong v_i F^\e_i$ as $F^\e_i$-module. Since $u_i\in
J_j$ for all $j>i$, one sees that $u_iu_j=u_i$ for all
$i<j$. This implies that $v_j v_i=v_j$ for all $j>i$. 
There are only two possibilities: either Case 1: $v_i=v_{i+1}$
for all $i\gg 0$, or Case 2: $v_{i}\neq v_{i+1}$ for infinitely
many $i$. In Case 1 we may assume that $v_i=v_{i+1}:=v$
for all $i$ by passing to a subsequence.
Then $F^\e_i=(1-v) F^\e_i\boplus v F^\e_i$ for all $i$.
Thus $F^\e=(1-v) F^\e\boplus v F^\e$. Since 
$$(1-v) F^\e=\bcup_{i}(1-v) F^\e_i=\bcup_i J_i=J,$$ 
$v F^\e\cong F$. This implies that $F$ is projective, a 
contradiction to Lemma \ref{xxlem7.2}. So Case 1 is impossible 
and that leaves us Case 2. By choosing a subsequence
of $\{F_i\}$ we may assume that $v_i\neq v_{i+1}$ for all
$i$. Let $e_1=u_1=1-v_1$ and $e_i=v_i-v_{i+1}$. Then 
$\{e_i\}$ is a set of nonzero orthogonal idempotents of
$F^\e$ and $J_n=u_n F^\e_n=\boplus_{n=1}^i e_i F^\e_i$ for all 
$n$. Since $J=\bcup_n J_n$, the assertion follows. 

(b) If $f: J\to F$ is 
a nonzero $F^\e$-homomorphism, then $f$ is surjective since $F$ is 
a field. Pick $b\in J_i\subset J$ such that $f(b)=1\in F$. Thus 
$f$ induces a nonzero $F_i^\e$-homomorphism from $bF_i^\e\to F_i$. 
But $F_i^\e$ is a direct sum $F_i\boplus J_i$ as rings, any homomorphism
from a submodule of $J_i$ to $F_i$ is zero. This is a contradiction,
hence $\Hom_{F_i^\e}(J,F)=0$. A similar argument shows that
$\Hom_{F_i^\e}(F,J)=0$.

(c) Applying $\Hom_{F^\e}(-,F)$ to the short exact sequence
\eqref{E7.2.1} we obtain an exact sequence
$$\to \Hom_{F^\e}(J,F)\to \Ext^1_{F^\e}(F,F)\to \Ext^1_{F^\e}(F^e,F)\to.$$
By part (b) the left end of the above sequence is zero and the right
end is zero since $F^\e$ is a free $F^\e$-module. Hence
$\Ext^1_{F^\e}(F,F)=0$.

If $\Hom_{F^\e}(F,F^e)\neq 0$, then let $F$ be the image of some 
nonzero map $F\to F^\e$. By part (b) $F\cap J=0$. Hence $F^\e=
F\boplus J$ because $F^\e/J\cong F$. This contradicts the fact
$\pdim_{F^\e} F=1$ in Lemma \ref{xxlem7.2}. Therefore 
$\Hom_{F^\e}(F,F^e)=0$. 

(d) Applying $\Hom_{F^\e}(F,-)$ to the short exact sequence
(\ref{E7.2.1}) we obtain an exact sequence
$$\to \Hom_{F^\e}(F,F)\to \Ext^1_{F^\e}(F,J)\to 
\Ext^1_{F^\e}(F, F^\e)\to.$$
Since $\Hom_{F^\e}(F,F)$ is 1-dimensional over $F$, it suffices to
show that $\Ext^1_{F^\e}(F,J)$ is infinite dimensional over $F$.
Recall that the short exact sequence \eqref{E7.2.1} is non-split
since $\pdim_{F^\e}F=1$. Hence \eqref{E7.2.1} represents a nonzero
element in $\Ext^1_{F^\e}(F,J)$, which we denote by $\psi$.

By part (a), $J=\boplus_{i=1}^{\infty} e_i F^\e$ where $\{e_i\}$ is 
an infinite set of nonzero orthogonal idempotents in $F^\e$.
Now let $\Phi$ be any infinite subset of ${\mathbb N}$ and
let $\Lambda=\boplus_{i\in \Phi} e_i F^\e$. We claim that 
$\Ext^1_{F^\e}(F,\Lambda)\neq 0$. Otherwise if 
$\Ext^1_{F^\e}(F,\Lambda)=0$, then 
$$\psi\in \Ext^1_{F^\e}(F,J)=\Ext^1_{F^\e}(F,J')$$
where $J'=\boplus_{i\not\in \Phi} \; e_i F^\e$. Hence
$\psi$ represents a non-split short sequence
$$0\to J'\to E\to F\to 0$$
such that
$$0\to J'\boplus \Lambda \to E\boplus \Lambda \to F\to 0$$
is equivalent to \eqref{E7.2.1}. Therefore $F^e\cong E\boplus 
\Lambda$; but this is impossible because $\Lambda$ is an infinite
direct sum. So we proved our claim that 
$\Ext^1_{F^\e}(F,\Lambda)\neq 0$. 

Next we decompose $\mathbb N$ into a disjoint union of 
infinitely many infinite subsets $\{\Phi_n\}_{n\in \mathbb N}$ and 
define $\Lambda_n
=\boplus_{i\in \Phi_n} e_i F^\e$. By the last paragraph,
$\Ext^1_{F^\e}(F, \Lambda_n)\neq 0$ for all $n$. Hence
the $F$-vector space dimension of $\Ext^1_{F^\e}(F, \boplus_{n=1}^p\Lambda_n)$ 
is at least $p$. Finally note that $\boplus_{n=1}^p\Lambda_n$
is a direct summand of $J$, therefore $\Ext^1_{F^\e}(F, J)$ is
infinite dimensional over $F$, as desired.
\end{proof}

\begin{proof}[Proof of Proposition \ref{xxprop7.1}]
(a) By Lemma \ref{xxlem7.2}, $F$ is smooth and $\shtr F\leq 1$.
By Lemma \ref{xxlem7.3}(d), $\injdim 
F^e>0$ and $\fhtr F>0$. The assertion follows from 
Lemma \ref{xxlem1.2}(a).

(b) Since $F$ is not finitely generated as a field, $F^e$
is not noetherian \cite[Proposition 1]{RSW}. 
We have seen that $F=\bcup_{i} G_i$
where $G_i$ ranges over all finite dimensional subfields of $F$.
For each $G_i$, we know that $\htr G_i=0$ [Example \ref{xxex1.9}(c)].
The assertion follows.

(c) Clearly $\tr F=0$. 

Let $U$ be 
a division algebra such that $F\otimes U$ is noetherian. We claim
that $F\otimes U$ is semisimple artinian. If the claim 
is proved, then $\injdim F\otimes U=\gldim F\otimes U=0$.
This implies that $\thtr F=0$. 

Now we prove the claim. Let $F=\bcup_i F_i$ where $\{F_i\}$ is 
an ascending chain of finite dimensional subfields of $F$. 
Then $F\otimes U=\bcup_i F_i\otimes U$. For each $i$, 
$F_i\otimes U$ is artinian since $F_i$ is finite dimensional.
Let $Z$ be the center of $U$. Then $F_i\otimes Z$ is a direct
sum of fields, say $\boplus_{t} G_t$, since $F_i$ is separable. 
Then $F_i\otimes U=\boplus_t G_t\otimes_Z U$. Since $Z$ is the
center of $U$, each $G_t\otimes_Z U$ is simple. Hence 
$F_i\otimes U$ is semisimple artinian. Write
\begin{equation}
\label{E7.3.1}
F_i\otimes U=\boplus_{t=1}^p \mathrm{M}_{n_t}(D^i_t)
\tag{E7.3.1}
\end{equation}
where $D^i_t$ are division rings. The Goldie rank of $F_i\otimes U$
is $\sum_{t=1}^p n_t$. Since $F\otimes U$ is faithfully flat over
$F_i\otimes U$, then Goldie rank of $F_i\otimes U$ is 
bounded by the Goldie rank of $F\otimes U$;
the latter is finite because $F\otimes U$ is noetherian. 
Therefore for $i\gg 0$, $F_i\otimes U$ has the same form of 
the decomposition \eqref{E7.3.1} with $D^i_t\subset D^{i+1}_t$ for all $i$.
Thus 
$$F\otimes U=\boplus_{t=1}^p \mathrm{M}_{n_t}(D_t)$$
where $D_t=\bcup_i D^i_t$. Since each $D^i_t$ is a division ring,
so it $D_t$. Therefore the claim is proved. 

(d) This follows from Lemma \ref{xxlem7.3}(d).
\end{proof}

\begin{example}
\label{xxex7.4} Let $F$ be the separable algebraic field extension 
of $k$ as in Proposition \ref{xxprop7.1}. Let $F'$ be another separable 
algebraic field extension of $k$ such that $F$ is a subfield of $F'$
and $F'$ is countably infinite dimensional over $F$. By Proposition
\ref{xxprop7.1}(a) $\htr F'=1$. So $\htr F'=\htr F$ but $\dim_F F'=\infty$.
\end{example}

\begin{example}
\label{xxex7.5} Let the base field $k$ be ${\mathbb C}(\{x_i\}_{i\in I})$
where $|I|>\aleph_m$ for all integers $m$. Let $F$ be the field
${\mathbb C}(\{x_i^{1/2}\}_{i\in I})$. Then $F$ is a separable algebraic 
field extension of $k$ such that $\dim_k F>\aleph_m$ for all integers $m$. 
For each $m$ there is a subfield $F_m\subset F$ such that
$\dim_k F_m=\aleph_m$. By \cite[Theorem 10]{Os}, $\pdim_{F^\e_m}F_m=
\gldim F^\e_m=m$. We claim that there is no smooth simple artinian 
ring $S$ such that $F\subset S$. If on contrary such $S$ exists, 
then $S^\e$ has finite global dimension and $S$ contains $F_m$ for 
all $m$. By \cite[Theorem 7.2.5]{MR}, $\gldim F^\e_m\leq \gldim S^\e$. 
Since $m$ is arbitrary, $\gldim S^\e=\infty$, a contradiction. 

It is unclear to us if $\htr F_m=m$ and if $F_m$ is homologically
uniform for all $m$. 
\end{example}

Next we are going to construct an algebra $A$ which is regular, but
not AS regular. 

Let $G$ be the nilpotent group generated by $a,b,c$ with relations
$ab=ba, ac=ca, bc=cba$. Then the group algebra $kG$ is isomorphic 
to $k[a^{\pm 1},b^{\pm 1}][c^{\pm 1},\sigma^{-1}]$ where $\sigma: 
a\mapsto a, b\mapsto ab$. The group algebra $kG$ is obviously 
$G$-graded. Let $A$ be the subalgebra of $kG$ generated by 
$x:=c,\; y:=ac,\; z:=bc, \; t:=abc$. Since $kG$ is a domain, so is 
$A$. If we set $\deg(c)=1$ and $\deg(a)=\deg(b)=0$, then $A$ is a 
connected $\mathbb N$-graded domain generated in degree $1$. 

The following proposition are due to a joint work of Paul Smith 
and the second author \cite{SmZ}. The authors thank Paul 
Smith for allowing them to use this unpublished result. 

\begin{proposition}\cite{SmZ}
\label{xxprop7.6} Let $A$ be the connected graded 
algebra constructed above.
\begin{enumerate}
\item
It is Koszul of global dimension $4$.
\item
It is a domain with $H_A(t)=(1-t)^{-4}$.
\item
It is not AS regular.
\item
It is neither left nor right noetherian.
\item
It has no non-trivial normal elements.
\end{enumerate}
\end{proposition}

\begin{remark}
\label{xxrem7.7} 
A part of the proof of Proposition \ref{xxprop7.6} 
was based on a long and tedious computation about 
the minimal free resolution of the trivial graded 
$A$-module $k$. It seems sensible to omit the proof here.
\end{remark}

We show now that Theorem \ref{xxthm0.2} is false without the
Artin-Schelter condition. Note that a domain of finite 
GK-dimension is an Ore domain. Hence the ring $A$ in Proposition
\ref{xxprop7.6} has a Goldie quotient ring. 

\begin{proposition}
\label{xxprop7.8}
Let $A$ be the algebra given in Proposition \ref{xxprop7.6}.
Let $Q(A)$ be the Goldie quotient ring of $A$. Then
$$\htr Q(A)=3< 4=\GKtr Q(A)=\GKtr A=\GKdim A.$$
\end{proposition}

\begin{proof} The division ring $Q(A)$ is isomorphic to
$k(a,b)(c;\sigma)$ where $\sigma$ maps $a\mapsto a,\; b\mapsto ba$.
So it is stratiform of length 3. By Proposition \ref{xxprop1.8},
$Q(A)$ is homologically uniform of $\htr$ $3$. The division algebra
$Q(A)$ is also the quotient division ring of the nilpotent group
$G=\langle a,b\rangle/(ab=ba,ac=ca,bc=cba)$. By a result of Lorenz
\cite[Theorem 2.2]{Lo}, $\GKtr Q(A)=4$. By 
\cite[Propositions 2.1 and 3.1(3)]{Zh1}, we have 
$\GKtr Q(A)\leq \GKtr A\leq \GKdim A=4$.
Therefore $\GKtr A=\GKdim A=4$. 
\end{proof}

Finally we post a question and make two remarks. 

\begin{question}
\label{xxque7.9} Let $S$ be a doubly noetherian simple artinian
algebra. Is then $S$ rigid and homologically uniform? 
\end{question}

\begin{remark}
\label{xxrem7.10} 
The work of Resco and Stafford \cite{Re1,Re2,Re3, St} suggests that 
one can define {\it Krull transcendence degree} of a simple artinian 
ring $S$ to be
$$\Ktr S=\Kdim S^\e$$
where $\Kdim$ denotes Krull dimension. Krull transcendence degree has 
nice properties similar to those listed in Proposition \ref{xxprop1.5}
(the proofs are also similar), we believe that this invariant deserves 
further study. 
\end{remark}

\begin{remark}
\label{xxrem7.11}
The division algebras in this paper are different from the 
free skew fields constructed by Cohen \cite{Co} and Schofield 
\cite{Sc2}. We expect that the free skew fields have infinite 
homological transcendence degree.
\end{remark}

\section*{Acknowledgments}
The authors thank Ken Goodearl, Tom Lenagan and Paul Smith for 
several useful discussions and valuable comments.

\providecommand{\bysame}{\leavevmode\hbox to3em{\hrulefill}\thinspace}
\providecommand{\MR}{\relax\ifhmode\unskip\space\fi MR }
\providecommand{\MRhref}[2]{%
  \href{http://www.ams.org/mathscinet-getitem?mr=#1}{#2} }
\providecommand{\href}[2]{#2}


\begin{thebibliography}{10}

\bibitem{Ar}
M. Artin,
\emph{Some problems on three-dimensional graded domains},
Representation theory and algebraic geometry 
(Waltham, MA, 1995), 1--19, London Math. Soc. Lecture 
Note Ser., 238, Cambridge Univ. Press, Cambridge, 1997.

\bibitem{ASZ}
M. Artin, L.W. Small and J.J. Zhang, 
\emph{Generic flatness for strongly Noetherian algebras},
J. Algebra \textbf{221} (1999), no. 2, 579--610.

\bibitem{ATV}
M. Artin, J. Tate and M. Van den Bergh,
\emph{Some algebras associated to automorphisms of 
elliptic curves},
The Grothendieck Festschrift, Vol. I, 33--85, 
Progr. Math., \textbf{86},
Birkh{\" a}user Boston, Boston, MA, 1990.

\bibitem{Co}
P.M. Cohn, 
\emph{Skew fields. Theory of general division rings}, 
Encyclopedia of Mathematics and its Applications, 57. 
Cambridge University Press, Cambridge, 1995. 

\bibitem{GK}
I.M. Gelfand and A.A. Kirillov,
\emph{Sur les corps li{\'e}s aux alg{\`e}bres 
enveloppantes des alg{\'e}bres de Lie}. 
(French) Inst. Hautes {\'E}tudes Sci. Publ. 
Math. No. 31 1966 5--19. 

\bibitem{GJ}
K.R. Goodearl and D.A. Jordan, 
\emph{Localizations of injective modules}, 
Proc. Edinburgh Math. Soc. (2) \textbf{28} 
(1985), no. 3, 289--299. 

\bibitem{GL}
K.R. Goodearl and T.H. Lenagan,
\emph{Catenarity in quantum algebras}, 
J. Pure Appl. Algebra 
\textbf{111} (1996), no. 1-3, 123--142. 

\bibitem{SGA6} 
A. Grothendieck et al., ``S\'eminare de G\'eom\'etrie 
Alg\'ebrique 6,'' LNM {\bf 225}, Springer 1971.

\bibitem{Lo}
M. Lorenz, 
\emph{On the transcendence degree of group 
algebras of nilpotent groups}, Glasgow Math. J. 
\textbf{25} (1984), no. 2, 167--174. 


\bibitem{MR}
J.C. McConnell and J.C. Robson, ``Noncommutative
        Noetherian Rings,'' Wiley, Chichester, 1987. 

\bibitem{Os}
B.L. Osofsky, 
\emph{Hochschild dimension of a separably 
generated field}, 
Proc. Amer. Math. Soc. {\bf 41} (1973), 24--30.  

\bibitem{Re1}
R. Resco, 
\emph{Transcendental division algebras and 
simple Noetherian rings},
Israel J. Math. \textbf{32} (1979), no. 2-3, 236--256. 

\bibitem{Re2}
R. Resco,
\emph{Krull dimension of Noetherian algebras 
and extensions of the base field}, 
Comm. Algebra \textbf{8} (1980), no. 2, 161--183. 
      
\bibitem{Re3}
R. Resco,
\emph{A dimension theorem for division rings}, 
Israel J. Math. 
\textbf{35} (1980), no. 3, 215--221. 

\bibitem{RSW}
R. Resco, L.W. Small and A.R. Wadsworth, 
\emph{Tensor products of division rings and 
finite generation of subfields}, 
Proc. Amer. Math. Soc. \textbf{77} (1979), 
no. 1, 7--10. 


\bibitem{RZ}
D. Rogalski and J.J. Zhang, 
\emph{Canonical maps to twisted rings}, 
arXiv:math.RA/0409405, preprint, 2004.


\bibitem{Ro}
A. Rosenberg,
\emph{Homological dimension and transcendence degree}, 
Comm. Algebra \textbf{10} (1982), no. 3, 329--338.


\bibitem{Rot}
J.J. Rotman, 
\emph{An introduction to homological algebra}, 
Pure and Applied Mathematics, 85. Academic Press, 
Inc. New York-London, 1979.  

\bibitem{Rou}
R. Rouquier,
\emph{Dimensions of triangulated categories}, 
arXiv:math.CT/0310134 v3, preprint, 2004.

\bibitem{Sc1}
A.H. Schofield,
\emph{Stratiform simple Artinian rings}, 
Proc. London Math. Soc. 
(3) \textbf{53} (1986), no. 2, 267--287.            
 
\bibitem{Sc2}
A.H. Schofield, 
\emph{Representation of rings over skew fields}, 
London Mathematical Society Lecture Note Series, 92. 
Cambridge University Press, Cambridge, 1985. 

\bibitem{Sc3}
A.H. Schofield, 
\emph{Questions on skew fields}, Methods in ring 
theory (Antwerp, 1983), 489--495, 
NATO Adv. Sci. Inst. Ser. C 
Math. Phys. Sci., 129, Reidel, Dordrecht, 1984.

\bibitem{SmZ}
P.S. Smith and J.J. Zhang, unpublished notes.

\bibitem{St}
J.T. Stafford, 
\emph{Dimensions of division rings}, 
Israel J. Math. \textbf{45} (1983), no. 1, 33--40. 

\bibitem{SZ}
J.T. Stafford and J.J. Zhang,
\emph{Algebras without Noetherian filtrations}, 
Proc. Amer. Math. Soc. \textbf{131} (2003), 
no. 5, 1329--1338.
 
\bibitem{SteZ}
D.R. Stephenson and J.J. Zhang, 
\emph{Growth of graded Noetherian rings}, 
Proc. Amer. Math. Soc. \textbf{125} (1997), 
no. 6, 1593--1605. 

\bibitem{Va}
P. V{\'a}mos, 
\emph{On the minimal prime ideal of a tensor 
product of two fields}, Math. Proc. Cambridge 
Philos. Soc. \textbf{84} (1978), no. 1, 25--35. 


\bibitem{VdB}
M. Van den Bergh, 
\emph{Existence theorems for dualizing complexes 
over non-commutative graded and filtered rings}, 
J. Algebra \textbf{195} (1997), no. 2, 662--679.

\bibitem{Ye1}
A. Yekutieli, 
\emph{Dualizing complexes over noncommutative 
graded algebras}, 
J. Algebra \textbf{153} (1992), no. 1, 41--84.

\bibitem{Ye2}
A. Yekutieli,
\emph{The residue complex of a noncommutative 
graded algebra},
J. Algebra \textbf{186} (1996), no. 2, 522--543. 

\bibitem{YZ1}
A. Yekutieli and J.J. Zhang,
\emph{Rings with Auslander dualizing complexes}, 
J. Algebra, \textbf{213} (1999), no. 1, 1--51.


\bibitem{YZ2}
A. Yekutieli and J.J. Zhang,
\emph{Residue complexes over noncommutative 
rings}, J. Algebra, \textbf{259} (2003), 
no. 2, 451--493. 

\bibitem{YZ3}
A. Yekutieli and J.J. Zhang, 
\emph{Dualizing complexes and perverse 
modules over Differential Algebras}, 
Compositio Math. \textbf{141} (2005), 
no. 3, 620--654.

\bibitem{YZ4}
A. Yekutieli and J.J. Zhang, 
\emph{Multiplicities of Indecomposable 
Injectives}, J. London Math. Soc. (2) 
\textbf{71} (2005), no. 1, 100--120.

\bibitem{YZ5}
A. Yekutieli and J.J. Zhang,
\emph{Dualizing Complexes and Tilting 
Complexes over Simple Rings},
J. Algebra, \textbf{256} (2002), no. 2, 556--567.

\bibitem{Zh1}
J.J. Zhang,
\emph{On Gelfand-Kirillov transcendence degree}, 
Trans. Amer. Math. Soc. \textbf{348} (1996), 
no. 7, 2867--2899. 

\bibitem{Zh2}
J.J. Zhang,
\emph{On lower transcendence degree}, Adv. Math. 
\textbf{139} (1998), no. 2, 157--193. 

\end{thebibliography}
\end{document}